\documentclass[12pt]{article}
\usepackage{amsfonts,a4,diss}
\usepackage{epsfig}
\ifx\onlineVersion\undefined%
\long\def\online#1#2{#2}%
\else%
\long\def\online#1#2{#1}%
\fi%
\newdimen\appletheight
\newdimen\appletwidth
\def\applet#1#2#3{%
\appletwidth=#2%
\appletheight=#3%
\vbox to\appletheight{
\vfill\hbox to\appletwidth{%
\vbox to \appletheight{}%
\hfill}}%
}

\def\paper{paper}
\def\thep{.}
\def\g{\gamma}
\def\F{\mathcal F}

\newcommand{\spann}{\mathop{\rm span}\nolimits}

\newcommand{\One}{\mathbb I}

\newcommand{\crt}{\mathop{\rm cr}\nolimits}

\renewcommand{\hat}{\widehat}

\begin{document}
\title{Discrete Hashimoto surfaces and a doubly discrete smoke ring
  flow} 
\author{Tim Hoffmann\\
} \maketitle
\begin{abstract}
  B\"acklund transformations for smooth and ``space discrete''
  Hashimoto surfaces are discussed and a geometric interpretation is
  given. It is shown that the complex curvature of a discrete space
  curve evolves with the discrete nonlinear Schr\"odinger equation
  (NLSE) of Ablowitz and Ladik, when the curve evolves with the
  Hashimoto or smoke ring flow. A doubly discrete Hashimoto flow is
  derived and it is shown, that in this case the complex curvature of
  the discrete curve obeys Ablovitz and Ladik's doubly discrete NLSE.
  Elastic curves (curves that evolve by rigid motion only under the
  Hashimoto flow) in the discrete and doubly discrete case are shown
  to be the same.

  \online{}{
    There is an online version of this paper, that can be viewed using
    any recent web browser that has JAVA support enabled. It includes
    two additional java applets. It can be found at
    \href{http://www-sfb288.math.tu-berlin.de/Publications/online/smokeringsOnline/index.html}%
    {
    http://www-sfb288.math.tu-berlin.de/Publications/online/smokeringsOnline/index.html}
  }
\end{abstract}
\label{sec:sr}
\section{Introduction}
\label{sr:sec:Introduction}
Many of the surfaces that can be described by integrable equations
have been discretized. Among them are surfaces of constant negative
Gaussian curvature, surfaces of constant mean curvature, minimal
surfaces, and affine spheres. This \paper\ continues the program by
adding Hashimoto surfaces to the list. These surfaces are obtained by
evolving a regular space curve $\g$ by the Hashimoto or {\em smoke
  ring flow}
\[ \dot\g = \g^\prime\times\g^{\prime\prime}\thep\]
As shown by Hashimoto \cite{HA72} this evolution is directly linked to
the famous nonlinear Schr\"odinger equation (NLSE)
\[ i \dot\Psi + \Psi^{\prime\prime} +\frac12\vert\Psi\vert^2\Psi = 0\thep\]
In \cite{AL76a} and \cite{AL76b} Ablowitz and Ladik gave a
differential-differen\-ce and a difference-difference discretization of
the NLSE. In \cite{HO99b}
the author shows\footnote{The equivalence for the differential-difference case
  appeared first in \cite{IS82}.}  that they correspond to a Hashimoto
flow on discrete curves (i.\ e.\ polygons) \cite{BS98,DS99} and a
doubly discrete Hashimoto flow respectively. This discrete evolution
is derived in section~\ref{sr:sec:DiscreteBaecklund} from a
discretization of the B\"acklund transformations for regular space
curves and Hashimoto surfaces.

In Section~\ref{sr:sec:smoothHashimoto} a short review of the smooth
Hashimoto flow and its connection to the isotropic Heisenberg magnet
model and the nonlinear Schr\"odinger equation is given. It is shown
that the solutions to the auxiliary problems of these integrable
equations serve as frames for the Hashimoto surfaces and a 
Sym formula is derived. In section~\ref{sr:sec:smoothAlg} the dressing
procedure or B\"acklund transformation is discussed and applied on
the vacuum. A geometric interpretation of this transformation as a
generalization of the Traktrix construction for a curve is given.

In Section~\ref{sr:sec:discreteHashimoto} the same program is carried
out for the Hashimoto flow on discrete curves.
Then in Section~\ref{sr:sec:ddHashimoto} special double B\"acklund
transformations (for discrete curves) are singled out to get a unique
evolution which serves as our doubly discrete Hashimoto flow.

Elastic curves (curves that evolve by rigid motion under the Hashimoto
flow) are discussed in all these cases. It turns out that discrete
elastic curves for the discrete and the doubly discrete Hashi\-moto flow
coincide.


Through this \paper\ we use a quaternionic description. Quaternions are
the algebra generated by $1,$ $\I,$ $\J,$ and $\K$ with the relations
$\I^2 = \J^2 = \K^2 = -1,$ $\I\J = K, \J\K = \I,$ and $\K\I = \J$.
Real and imaginary part of a quaternion are defined in an obvious
manner: If $q = \alpha + \beta\I + \gamma\J + \delta\K$ we set $\Re(q)
= \alpha$ and $\Im(q) = \beta\I + \gamma\J + \delta\K$. Note that
unlike in the complex case the imaginary part is not a real number.
We identify the 3-dimensional euclidian space with the imaginary
quaternions i.~e.\ the span of $\I,$ $\J,$ and $\K$. Then for two
imaginary quaternions $q,r$ the following formula holds:
\[ qr = -\<q,r> + q\times r\] with $\<\cdot,\cdot>$ and
$\cdot\times\cdot$ denoting the usual scalar and cross products of
vectors in 3-space. A rotation of an imaginary quaternion around the
axis $r, \vert r\vert = 1$ with angle $\phi$ can be written as
conjugation with the unit length quaternion $(\cos\frac\phi2 +
\sin\frac\phi2r)$. 

Especially when dealing with the Lax representations of the various
equations it will be convenient to identify the quaternions with
complex 2 by 2 matrices:
\[ \I = i\sigma_3 = \quadmatrix{i}{0}{0}{-i}\quad \J = i\sigma_1 =
\quadmatrix{0}{i}{i}{0}\]
\[\K = -i\sigma_2 = \quadmatrix{0}{-1}{1}{0}\thep\]

\section[Hashimoto flow, Heisenberg flow, and the NLSE]{The Hashimoto flow, the Heisenberg flow and the nonlinear
  Schr\"odinger equation}
\label{sr:sec:smoothHashimoto}
Let $\g:\R\to \R^3 = \Im\H$ be an arclength parametrized regular curve
and $\F:\R\to\H^*$ be a parallel frame for it, i.\ e.\
\begin{eqnarray}
\F^{-1}\I\F&=&\g^\prime = \g_x\\
(\F^{-1}\J\F)^\prime&\|& \g^\prime\thep\label{sr:parallelFrameCond}
\end{eqnarray}
The second equation says that $\F^{-1}\J\F$ is a parallel section in
the normal bundle of $\g$. which justifies the name.  Moreover let $A
= \F^\prime\F^{-1}$ be the logarithmic derivative of $\F$.  Equation
(\ref{sr:parallelFrameCond}) gives, that $A$ must lie in the
$\J$-$\K$-plane and thus can be written as
\begin{equation}
        A = -\frac\Psi2\K
        \label{sr:complexCurvature}
\end{equation}
with  $\Psi\in\spann(1,\I)\cong\C$.
\begin{definition}
We call $\Psi$ the complex curvature of $\g$.
\end{definition}

Now let us evolve $\g$ with the following flow:
\begin{equation}
        \dot\g = \g^\prime\times\g^{\prime\prime} = \g^\prime\g^{\prime\prime}\thep
        \label{sr:eq:smokreRingFlow}
\end{equation}
Here $\dot\g$ denotes the derivative in time. This is an evolution in
binormal direction with velocity equal to the (real) curvature. It is
known as the {\em Hashimoto}\/ or {\em smoke ring flow}. Hashimoto was
the first to show, that under this flow the complex curvature $\Psi$
of $\g$ solves the nonlinear Schr\"odinger equation (NLSE) \cite{HA72}
\begin{equation}
        \I\dot\Psi + \Psi^{\prime\prime} +\frac12\vert\Psi\vert^2\Psi = 0\thep
        \label{sr:nlse}
\end{equation}
or written for $A$:
\begin{equation}
  \label{sr:nlseForA}
  \I\dot A + A^{\prime\prime} = 2A^3\thep
\end{equation}
\begin{definition}
  The surfaces $\g(x,t)$ wiped out by the flow
   given in equation (\ref{sr:eq:smokreRingFlow}) are called Hashimoto
  surfaces.
\end{definition}
Equation (\ref{sr:nlse}) arises as the zero curvature condition $\hat L_t -
\hat M_x + [\hat L,\hat M]=0$ of the
system
\begin{equation}
  \label{sr:eq:smoothAux}
\begin{array}{rcl}
  \hat\F_x(\mu) & = & \hat L(\mu)\hat\F(\mu)\\
  \hat\F_t(\mu) & = & \hat M(\mu)\hat\F(\mu)
\end{array}
\end{equation}
with 
\begin{equation}
  \label{sr:eq:nlseLaxPair}
\begin{array}{rcl}
  \hat L(\mu) & = & \mu\I -\frac\Psi2\K\\
  \hat M(\mu) & = & \frac{\vert\Psi\vert^2}4\I + \frac{\Psi_x}2\J -
  2 \mu\hat L(\mu)\thep
\end{array}
\end{equation}
To make the connection to the description with the parallel frame $\F$ 
we add torsion to the curve $\g$ by setting 
\[ A(\mu) = e^{-2 \mu x\I}\Psi\K\thep\]
This gives rise to a family of curves $\g(\mu)$ the so-called {\em associated
family}\/ of $\g$. Now one can gauge the corresponding parallel frame
$\F(\mu)$ with $e^{\mu x\I}$ and get 
\[ (e^{\mu x\I}\F(\mu))_x = ((e^{\mu x\I})_x e^{-\mu x\I} +e^{\mu
  x\I}A(\mu)e^{-\mu x\I})e^{\mu x\I}\F(\mu) = L(\mu)e^{\mu
  x\I}\F(\mu)\] with $L(\mu)$ as in (\ref{sr:eq:smoothAux}). So above
$\hat\F(\mu)= e^{\mu x\I}\F(\mu) $ is for each $t_0$ a frame for the
  curve $\g(x,t_0)$.
\begin{theorem}{\rm \bf (Sym formula)}
  Let $\Psi(x,t)$ be a solution of the NLSE (equation
  (\ref{sr:nlse})). Then up to an euclidian motion the corresponding
  Hashimoto surface $\g(x,t)$ can be obtained by
  \begin{equation}
    \label{sr:eq:smoothSym}
    \g(x,t) = \hat\F^{-1}\hat\F_\lambda\vert_{\lambda=0}
  \end{equation}
  where $\hat\F$ is a solution to (\ref{sr:eq:smoothAux}).
\end{theorem}
\begin{proof}
  Obviously $\hat\F\vert_{\lambda=0}(x,t_0)$ is a parallel frame for each $\g(x,t_0)$.
  So writing $\hat\F(x,t_0)\vert_{\lambda=0} =: \F(x),$
  one easily computes $(\hat\F^{-1}\hat\F_\lambda\vert_{\lambda=0})_x
  = \F^{-1}\I\F=\g_x$ and
  $(\hat\F^{-1}\hat\F_\lambda\vert_{\lambda=0})_y = \F^{-1}\Psi\K\F$.
  But $\g_t = \g_x\g_{xx} = \F^{-1}\Psi\K\F$.
\end{proof}

If one differentiates equation (\ref{sr:eq:smokreRingFlow}) with respect
to $x$ one gets the so-called isotropic Heisenberg magnet model (IHM):
\begin{equation}
  \label{sr:eq:heisenberg}
  \dot S = S\times S^{\prime\prime} = S\times S_{xx}
\end{equation}
with $S = \g^\prime$. This equation arises as zero curvature condition
$U_t -V_x + [U,V] = 0$ with matrices
\begin{equation}
  \label{sr:eq:smoothIHMLaxPair}
  \begin{array}{rcl}
    U(\lambda) &=& \lambda S\\
    V(\lambda) &=& -2 \lambda^2 S -\lambda S^\prime S
  \end{array}
\end{equation}
In fact if $G$ is a solution to
\begin{equation}
  \label{sr:eq:smoothIHMLinearProb}
  \begin{array}{rcl}
    G_x = U(\lambda) G\\
    G_t = V(\lambda) G
  \end{array}
\end{equation}
it can be viewed as a frame for the Hashimoto surface too and one has a
similar Sym formula:
\begin{equation}
  \label{sr:eq:heisenbergSymFormula}
  \g(x,t) =  G^{-1}G_\lambda\vert_{\lambda=0}
\end{equation}
The system (\ref{sr:eq:smoothIHMLinearProb}) is known to be gauge
equivalent to (\ref{sr:eq:smoothAux}) \cite{FT86}.

\subsection{Elastic curves}
\label{sr:sec:ElasticCurves}

The stationary solutions of the NLSE (i.\ e.\ the curves that evolve
by rigid motion under the Hashimoto flow) are known to be the {\em
  elastic curves}\/ \cite{BS98}. They are the critical points of the
functional
\[E(\g) = \int \kappa^2\]
with $\kappa = \vert\Psi\vert$ the curvature of $\g$.  The fact that
they evolve by rigid motion under the Hashimoto flow can be used to
give a characterization by their complex curvature $\Psi$ only: When
the curve evolves by rigid motion $\Psi$ may get a phase factor only.
Thus $\dot\Psi = c \I \Psi$. Inserted into equation~(\ref{sr:nlse})
this gives
\begin{equation}
  \label{sr:eq:PsiforElastica}
  \Psi^{\prime\prime} = (c-\frac12\vert\Psi\vert^2)\Psi\thep
\end{equation}

\subsection{B\"acklund transformations for smooth space curves and
  Hashimoto surfaces} Now we want to describe the dressing procedure
or B\"acklund transformation for the IHM model and the Hashimoto
surfaces. This is a method to generate new solutions of our equations
from a given one in a purely algebraic way. Afterwards we give some
geometric interpretation for this transformation.

\subsubsection{Algebraic description of the B\"acklund transformation}
\label{sr:sec:smoothAlg}
\begin{theorem}
\label{sr:thm:smoothDressing}
  Let $G$ be a solution to 
  equations (\ref{sr:eq:smoothIHMLinearProb})
  with $U$ and $V$ as in (\ref{sr:eq:smoothIHMLaxPair}) (i.\ e.\ $U(1)$
  solves the $IHM$ model).
  Choose $\lambda_0,s_0\in\C$. Then $ \tilde G(\lambda) :=
  B(\lambda)G(\lambda)$ with $B(\lambda) = (\One + \lambda \rho),
  \rho\in\H$ defined by the conditions that $\lambda_0,\bar\lambda_0$
  are the zeroes of $\det(B(\lambda))$ and 
  \begin{equation}
    \label{sr:eq:BTcondition}
    \tilde G(\lambda_0){s_0\choose 1} = 0 \quad\hbox{\rm and}\quad \tilde
  G(\bar\lambda_0){1\choose-\bar s_0} = 0
  \end{equation}
  solves a system of the same type. In particular $\tilde U(1) =
  \tilde G_x(1)\tilde G^{-1}(1)$ solves again the Heisenberg magnet model
  (\ref{sr:eq:heisenberg}).
\end{theorem}
\begin{proof}
  We define $\tilde U(\lambda) = \tilde G_x\tilde G^{-1}$ and $\tilde
  V(\lambda) = \tilde G_t\tilde G^{-1}$. Equation
  (\ref{sr:eq:BTcondition}) ensures that $\tilde U(\lambda)$ and
  $\tilde V(\lambda)$ are smooth at $\lambda_0$ and $\bar\lambda_0$.
  Using $\tilde U(\lambda) = B_x(\lambda)B^{-1}(\lambda) +
  B(\lambda)U(\lambda)B^{-1}(\lambda)$ this in turn implies that
  $\tilde U(\lambda)$ has the form $\tilde U(\lambda) = \lambda \tilde
  S$ for some $\tilde S$.

  Since the zeroes of $\det(B(\lambda))$ are fixed we know that
  $r := \Re(\rho)$ and $l := \vert\Im(\rho)\vert$ are constant. We write
  $\rho = r + v$.

  One gets $\tilde S = S + v_x$ and 
  \begin{equation}
    \label{sr:eq:vx}
    v_x = \frac{2rl}{r^2 + l^2}\frac{v\times S}{l} + \frac{2l^2}{r^2 +
    l^2}\frac{\<v,S>}{l^2}v - \frac{2l^2}{r^2 + l^2}S\thep
  \end{equation}
  This can be used to show $\vert \tilde S\vert = 1$.
  
  Again equation (\ref{sr:eq:BTcondition}) ensures that $\tilde
  V(\lambda) = \lambda^2 X + \lambda Y$ for some $X$ and $Y$. But then
  the integrability condition $\tilde U_t - \tilde V_x + [\tilde U,
  \tilde V]$ gives up to a factor $c$ and possible constant real
  parts $x$ and $y$ that $X$ and $Y$ are fixed to be $X =x + c\tilde
  S_x\tilde S + d\tilde S$ and $Y = y + 2c\tilde S$. The additional
  term $d\tilde S$ in $X$ corresponds to the (trivial) tangential flow
  which always can be added. The form $\tilde V(\lambda) =
  B_t(\lambda)B^{-1}(\lambda) + B(\lambda)V(\lambda)B^{-1}(\lambda)$
  gives $c=-1$ and $d=0$. Thus one ends up with $\tilde V(\lambda) =
  -2\lambda^2 \tilde S -\lambda \tilde S_x\tilde S$.
\end{proof}
So we get a four parameter family ($\lambda_0$ and $s_0$ give two
real parameters each) of transformations for our curve $\g$ that are
compatible with the Hashimoto flow. They correspond to the four
parameter family of B\"acklund transformations of the
NLSE.

\begin{figure}[tbp]
  \begin{center}
    \epsfxsize=0.55\hsize\epsfbox{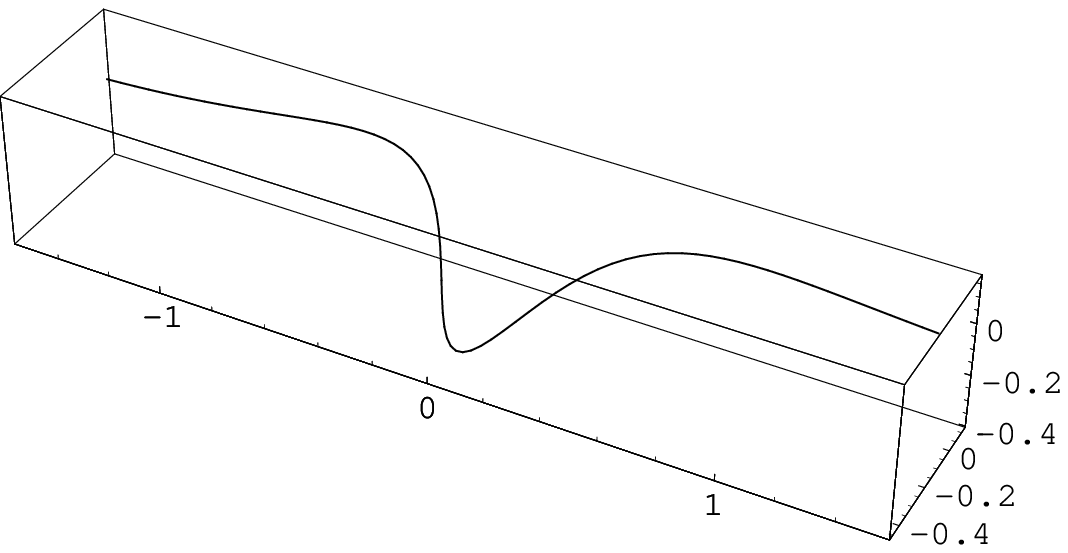}\epsfxsize=0.55\hsize\kern-0.1\hsize\raise7pt\hbox{\epsfbox{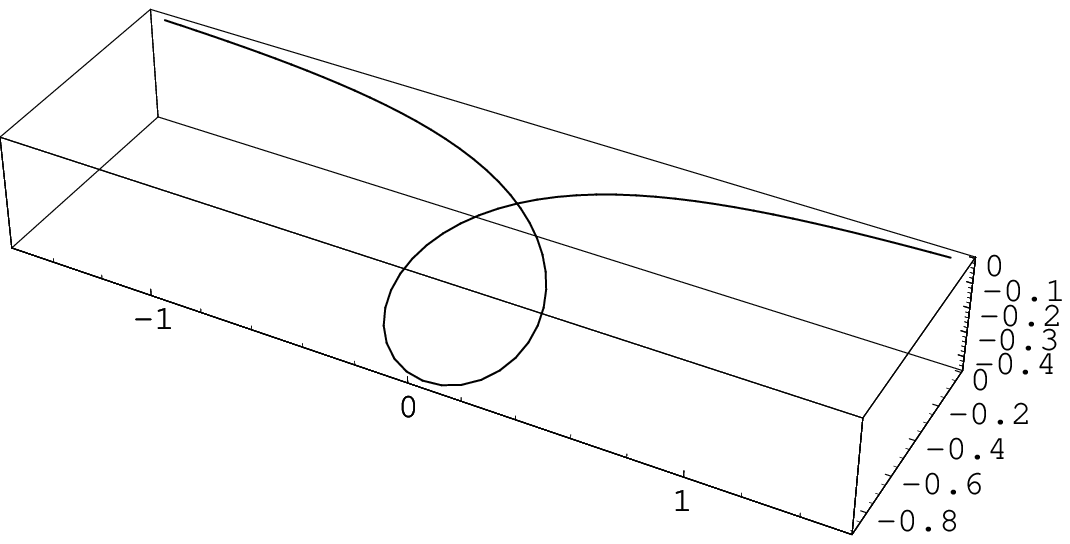}}

    \epsfxsize=0.55\hsize\epsfbox{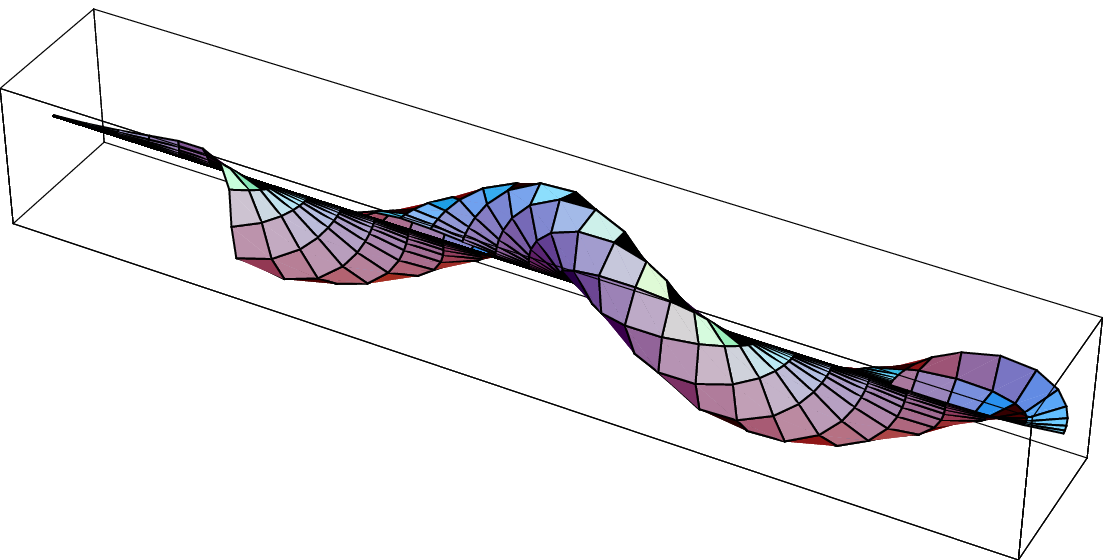}\epsfxsize=0.55\hsize\kern-0.1\hsize\raise7pt\hbox{\epsfbox{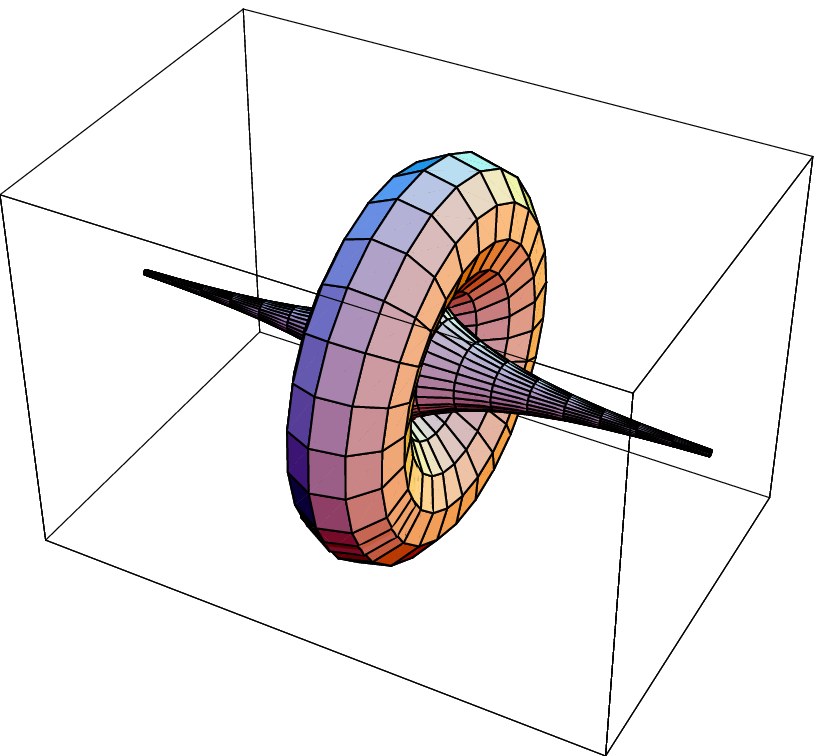}}
    \caption{Two dressed straight lines and the corresponding
      Hashimoto surfaces}
    \label{sr:fig:dressedLine}
  \end{center}
\end{figure}
\begin{example}
  Let us do this procedure in the easiest case:
  We choose $S \equiv \I$ (or $\g(x,t) = x \I$) which gives 
  \[G(\lambda) =\exp((\lambda x - 2\lambda^2 t) \I) = \quadmatrix{e^{i(\lambda
      x- 2\lambda^2 t)}}{0}{0}{e^{-i(\lambda x- 2\lambda^2 t)}}\thep\]
  After choosing $\lambda_0$ and $s_0$ and writing $\rho =
  \quadmatrix{a}{b}{-\bar b}{\bar a}$ one gets with equation
  (\ref{sr:eq:BTcondition}) 
  \begin{equation}
    \label{sr:eq:example1}
    \begin{array}{rcl}
      -e^{i (\lambda_0 x - 2 \lambda_0^2 t)} & = & \lambda_0(e^{i (\lambda_0 x - 2 \lambda_0^2 t)} a +
      s_0 e^{-i (\lambda_0 x - 2 \lambda_0^2 t)} b\\
      s_0 e^{-i (\lambda_0 x - 2 \lambda_0^2 t)} & = & \lambda_0(e^{i (\lambda_0 x - 2 \lambda_0^2 t)}\bar b
      - s_0 e^{-i (\lambda_0 x - 2 \lambda_0^2 t)} \bar a).
    \end{array}
  \end{equation}
  These equations can be solved for $a$ and $b:$
  \begin{equation}
    \label{sr:eq:example2}
    \begin{array}{rcl}
      a & = & -\frac{\frac1{\lambda_0} + \frac{s_0\bar
          s_0}{\bar\lambda_0}e^{-2i(\lambda_0-\bar\lambda_0)x+4i(\lambda_0^2-\bar\lambda_0^2)t}}%
      {1 + s_0\bar s_0 e^{-2i(\lambda_0-\bar\lambda_0)x+4i(\lambda_0^2-\bar\lambda_0^2)t}}\\
      b & = & \bar s_0 e^{2i\bar\lambda_0 x-4i\bar\lambda_0^2t}\frac{\frac1{\bar
          \lambda_0} - \frac1{\lambda_0}}%
      {1 + s_0\bar s_0 e^{-2i(\lambda_0-\bar\lambda_0)x+4i(\lambda_0^2-\bar\lambda_0^2)t}}
    \end{array}
  \end{equation}
  Using the Sym formula (\ref{sr:eq:heisenbergSymFormula}) one can
  immediately write the formula for the resulting Hashimoto surface
  $\tilde \g$:
  \[ \tilde \g = \Im(\rho) +\g = \quadmatrix{\Im(a) + i x}{b}{-\bar
    b}{-\Im(a) - i x}.\]
  The need for taking the imaginary part is due to the fact that we
  did not normalize $B(\lambda)$ to $\det(B(\lambda))=1$.
  
  If one wants to have the result in a plane $\arg b$ should be
  constant. This can be achieved by choosing $\lambda\in i\R$.
  Figure~\ref{sr:fig:dressedLine} shows the result for $s_0 = 0.5 + i$
  and $\lambda_0 = 1 - i$ and $\lambda_0 = - i$ respectively.
\end{example}

Of course one can iterate the dressing procedure to get new curves (or
surfaces) and it is a natural question how many one can get. This
leads immediately to the Bianchi permutability theorem
\begin{theorem}{\rm\bf (Bianchi permutability)}
  \label{sr:thm:smoothBianchi}
  Let $\tilde\g$ and $\hat\g$ be two B\"ack\-lund transforms of
  $\g$. Then there is a unique Hashimoto surface $\hat{\tilde\g}$ that
  is B\"acklund transform of $\tilde\g$ and $\hat\g$.
\end{theorem}
\begin{proof}
  Let $G,\hat G,$ and $\tilde G$ be the solutions to
  (\ref{sr:eq:smoothIHMLinearProb}) corresponding to $\g,\hat\g,$ and
  $\tilde\g$. One has $\hat G = \hat B G$ and $\tilde G = \tilde B G$
  with $\hat B = \One + \lambda\hat\rho$ and $\tilde B = \One +
  \lambda\tilde\rho$. The ansatz $\tilde{\hat B} \hat G = \hat{\tilde
    B}\tilde G$ leads to the compatability condition
  $\tilde{\hat B}\hat B = \hat{\tilde B}\tilde B$ or
  \begin{equation}
    \label{sr:eq:smoothPermutability}
    (\One + \lambda\tilde{\hat\rho})(\One +\lambda\hat\rho) = (\One
    +\lambda\hat{\tilde\rho})(\One + \lambda\tilde\rho)
  \end{equation}
  which gives:
  \begin{equation}
    \label{sr:eq:smoothPermEvol}
    \begin{array}{rcl}
      \tilde{\hat\rho} & = & (\hat\rho - \tilde\rho)\; \tilde\rho\; (\hat\rho
      - \tilde\rho)^{-1}\\
      \hat{\tilde\rho} & = & (\hat\rho - \tilde\rho)\; \hat\rho\; (\hat\rho
      - \tilde\rho)^{-1}\thep
    \end{array}
  \end{equation}
  Thus $\tilde{\hat B}$ and $\hat{\tilde B}$ are completely
  determined. To show that they give dressed solutions we note that
  since $\det\tilde{\hat B}\det\hat B = \det\hat{\tilde B}\det\tilde
  B$ the zeroes of $\det\tilde{\hat B}$ are the same as the ones of
  $\det\tilde B$ (and the ones of $\det\hat{\tilde B}$ coincide with
  those of $\det\hat B$). Therefore they do not depend on $x$ and
  $t$. Moreover at these points the kernel of $\tilde{\hat B}\hat G$ 
  coincides with the one of $\tilde G$. Thus it does not depend
  on $x$ or $t$ either. Now theorem~\ref{sr:thm:smoothDressing} gives
  the desired result.
\end{proof}
\subsubsection{Geometry of the B\"acklund transformation}

As before let $\g:I\to\R^3 = \Im\H$ be an arclength parametrized
regular curve. Moreover let $v:I\to\R^3 = \Im\H,$ $\vert v\vert = l$
be a solution to the following system:
\begin{equation}
  \label{sr:eq:traktrix}
  \begin{array}{rcl}
    \hat\g &=& \g + \frac12v\\
    \hat\g^\prime &\|& v\thep
  \end{array}
\end{equation}
Then $\hat\g$ is called a Traktrix of $\g$.
The forthcoming definition in this section is motivated by the
following observation: 
If we set $\tilde\g = \g + v$ it is again an arclength parametrized
curve and $\hat\g $ is a Traktrix of $\tilde\g$  too.
One can generalize this in the following way:
%
%
%
\begin{lemma}\label{sr:contConstLengthLemma}
  Let $v:I\to\Im\H$ be a vector field along $\g$ of constant length $l$ satisfying
  \begin{equation}
    \label{sr:eq:vx2}
    v^\prime = 2\sqrt{b - b^2}\; \frac{v\times\g^\prime}l + 2b 
    \frac{<v,\g^\prime>}{l^2} v - 2b\g^\prime
  \end{equation}
  with $0\leq b\leq 1$.
  Then $\tilde\g = \g + v$ is arclength parametrized.
\end{lemma}
\begin{proof}
  Obviously the above transformation coincides with the dressing
  described in the last section with $b = \frac{l^2}{r^2 + l^2}$ in
  formula (\ref{sr:eq:vx}). This proves the lemma.
\end{proof}

So $\Im(\rho)$ from theorem~\ref{sr:thm:smoothDressing} is nothing but
the difference vector between the original curve and the
B\"acklund transform.
Note that in the case $b = 1$ one gets the above Traktrix construction, 
that is for $\hat\g = \g + v$ holds $\hat\g^\prime\parallel v$. This
motivates the following 
\begin{definition}\label{sr:def:twistedTraktrix}
  The curve $\hat\g = \g + \frac12v$ with v as in
  lemma~\ref{sr:contConstLengthLemma} is called a {\em twisted
    Traktrix}\/ of the curve $\g$ and $\tilde\g = \g + v$ is called a
  {\em B\"acklund transform}\/ of $\g$.
\end{definition}

Moreover equation (\ref{sr:eq:vx2}) gives that $v'\perp v$ and
therefore $\vert v\vert \equiv const$. Since $v = \tilde\g - \g$ we
see that the B\"acklund transform is in constant distance to the
original curve.

\section[Discr. Hashimoto flow, Heisenberg flow, and dNLSE]{The Hashimoto flow, the Heisenberg flow, and
  the nonlinear Schr\"odinger equation in the discrete case}
\label{sr:sec:discreteHashimoto}
In this section we give a short review on the discretization (in
space) of the Hashimoto flow, the  isotropic Heisenberg magnetic
model, and the nonlinear Schr\"odinger equation. For more details on
this topic see \cite{FT86,BS98,DS99} and \cite{HO99b}. 

We call a map $\g:\Z\to\Im\H$ a discrete regular curve if any two
successive points do not coincide. It will be called arclength
parametrized curve, if $\vert \g_{n+1} - \g_n\vert = 1$ for all
$n\in\Z$. We will use the notation $S_n :=\g_{n+1}-\g_n$.  The
binormals of the discrete curve can be defined as $\frac{S_{n}\times
  S_{n-1}}{\vert S_{n}\times S_{n-1}\vert}$.

There is a natural discrete analog of a parallel frame:
\begin{definition}
A {\em discrete parallel frame}\/ is a map $\F:\Z\to\H^*$ with $\vert
\F_k\vert = 1$ satisfying
\begin{eqnarray}
        S_n & = & \F_n^{-1} \I\F_n
        \label{sr:discreteFrame1}  \\
        \Im\left((\F_{n+1}^{-1}\J\F_{n+1})(\F_{n}^{-1}\J\F_{n})\right) & \| 
        &\Im\left(S_{n+1}S_n\right)\thep
        \label{sr:discreteFrame2}
\end{eqnarray}  
\end{definition}
Again we set $\F_{n+1}  =  A_n\F_n$ and in complete analogy to the
continuous case eqn (\ref{sr:discreteFrame2}) gives the following
form for $A:$
\[A = \cos\frac{\phi_n}2 - \sin\frac{\phi_n}2\exp\left(\I\sum_{k=0}^n
  \tau_k\right) \K\] with $\phi_n = \angle\left(S_n,S_{n+1}\right)$
the folding angles and $\tau_n$ the angles between successive
binormals.  If we drop the condition that $\F$ should be of unit
length we can renormalize $A_n$ to be $1 -
\tan\frac{\phi_n}2\exp\left(\I\sum_{k=0}^n \tau_k\right) \K =: 1 -
\Psi_n\K$ with $\Psi_n\in\spann(1,\I)\cong\C$ and $\vert\Psi_n\vert =
\kappa_n$ the discrete (real) curvature.
\begin{definition}
  We call $\Psi$ the complex curvature\footnote{It would be
    more reasonable to define $A = 1 - \frac{\Psi_n}2\K$. which
    implies $\kappa_n = 2\tan\frac{\phi_n}2$ but notational simplicity
    makes the given definition more convenient.}
  of the discrete curve $\g$.
\end{definition}

Discretizations of the Hashimoto flow (\ref{sr:eq:smokreRingFlow})
(i.\ e.\ a Hashimoto flow for a discrete arclength parametrized curve) 
and the isotropic Heisenberg model (eqn (\ref{sr:eq:heisenberg})) are
well known \cite{FT86} (see also \cite{BS98} for a good discussion of
the topic). In particular a discrete version of
(\ref{sr:eq:smokreRingFlow}) is given by:
\begin{equation}
  \label{sr:eq:discreteSmokeRingFlow}
  \dot\g_k = 2\frac{S_{k}\times S_{k-1}}{1 + \<S_{k},S_{k-1}>}
\end{equation}
which implies for a discretization of (\ref{sr:eq:heisenberg})
\begin{equation}
  \label{sr:eq:discreteHeisenberg}
  \dot S_k = 2\frac{S_{k+1}\times S_k}{1 +\<S_{k+1},S_k>} -
  2\frac{S_{k}\times S_{k-1}}{1 +\<S_{k},S_{k-1}>}
\end{equation}
Let us state the zero curvature representation for this equation too:
Equation (\ref{sr:eq:discreteHeisenberg}) is the compatibility
condition of  $\dot U_k = V_{k+1}U_k - U_kV_k$ with
\begin{equation}
  \label{sr:eq:discreteIHMLaxPair}
  \begin{array}{rcl}
    U_k &=& \One + \lambda S_k\\[0.2cm]
    V_k &=& -\frac1{1+\lambda^2}\left(2 \lambda^2 \frac{S_{k} +
    S_{k-1}}{1 + \<S_{k},S_{k-1}>} + 2 \lambda \frac{S_{k}\times
    S_{k-1}}{1 +\<S_{k},S_{k-1}>} \right)
  \end{array}
\end{equation}
The solution to the auxiliary problem
\begin{equation}
  \label{sr:eq:dIHMAuxProb}
  \begin{array}{rcl}
    G_{k+1} & = & U_k(\lambda)G_k\\
    \dot G_k & = & V_k(\lambda)G_k
  \end{array}
\end{equation}
can be viewed as the frame to a discrete Hashimoto surface $\g_k(t)$ and
one has the same Sym formula as in the continuous case:
\begin{theorem}
  Given a solution $G$ to the system (\ref{sr:eq:dIHMAuxProb}) the
  corresponding discrete Hashimoto surface can be obtained up to an
  euclidian motion by 
  \begin{equation}
    \label{sr:eq:dIHMSym}
    \g_k(t) = (G_k^{-1}\;\frac\partial{\partial
    \lambda}G_k)\vert_{\lambda=0}. 
  \end{equation}
\end{theorem}
\begin{proof}
  One has $G_k^{-1}\;\frac\partial{\partial
    \lambda}G_k\vert_{\lambda=0} = \sum_{i=0}^{k-1} S_i =
    \g_k$ for fixed time $t_0$ and 
    \[(G_k^{-1}\frac\partial{\partial
    \lambda}G_k\vert_{\lambda=0})_t = (\frac\partial{\partial
    \lambda}V_k(\lambda)\vert_{\lambda=0}) = 2\frac{S_{k}\times
    S_{k-1}}{1 + \<S_{k},S_{k-1}>}\thep\]
\end{proof}

To complete the analogy to the smooth case we give a discretization of
the NLSE that can be found in \cite{AL76a} (see also
\cite{FT86,SU97}):
  \begin{equation}
    \label{sr:eq:dNLSE}
    -i\dot\Psi_k = \Psi_{k+1} - 2 \Psi_k + \Psi_{k-1} +\vert\Psi_k\vert^2
  (\Psi_{k+1} + \Psi_{k-1})\thep
  \end{equation}
\begin{theorem}
  \label{sr:thm:equivalenceDIHMdNLSE}
  Let $\g$ be a discrete arclength parametrized curve. If $\g$ evolves
  with the discrete Hashimoto flow (\ref{sr:eq:discreteSmokeRingFlow})
  then its complex curvature $\Psi$ evolves with the discrete
  nonlinear Schr\"odinger equation (\ref{sr:eq:dNLSE})
\end{theorem}

A proof of this theorem can be found in \cite{IS82} and
\cite{HO99b}.  There is another famous discretization of
the NLSE in literature that is related to the dIHM \cite{IK81,FT86}.
Again in \cite{HO99b} it is shown that it is in fact gauge
equivalent to the above cited which turns out to be more natural from
a geometric point of view.

\subsection{Discrete elastic curves}
\label{sr:sec:discreteElastic}

As mentioned in Section~\ref{sr:sec:ElasticCurves} the stationary
solutions of the NLSE (i.\ e.\ the curves that evolve by rigid motion
under the Hashimoto flow) are known to be the {\em elastic curves.}\/
They have a natural discretization using this property:

\begin{definition}
  A discrete  elastic curve is a curve $\g$ for which the evolution of
  $\g_n$ under the Hashimoto flow (\ref{sr:eq:discreteSmokeRingFlow}) is a
  rigid motion which means that its tangents evolve under the discrete
  isotropic Heisenberg model
  (\ref{sr:eq:discreteHeisenberg}) by rigid rotation.
\end{definition}
In \cite{BS98} Bobenko and Suris showed the equivalence of this
definition to a variational description.

The fact that  (\ref{sr:eq:discreteHeisenberg}) has to be a rigid rotation
means that the left hand side must be $S_n\times p$ with a unit
imaginary quaternion $p$. We will now give a description of elastic
curves by their complex curvature function only:

\begin{theorem}
  The complex curvature $\Psi_n$ of a discrete elastic curve $\g_n$
  satisfies the following difference equation:
  \begin{equation}
    \label{sr:eq:elasticCurvature}
    {\mathcal C} \frac{\Psi_n}{1 + \vert\Psi_n\vert^2} = \Psi_{n+1} + \Psi_{n-1}
  \end{equation}
  for some real constant ${\mathcal C}$.
\end{theorem}
Equation (\ref{sr:eq:elasticCurvature}) is a special case of a
discrete-time Garnier system (see \cite{SU94}).

\begin{proof}
  One can proof the theorem by direct calculations or using the
  equivalence of the dIHM model and the dNLSE stated in
  theorem~\ref{sr:thm:equivalenceDIHMdNLSE}. If the curve $\g$ evolves
  by rigid motion its complex curvature may vary by a phase factor
  only: $\Psi(x,t) = e^{i\lambda(t)}\Psi(x,t_0)$ or $\dot\Psi =
  i\dot\lambda\Psi$. Plugging this in eqn (\ref{sr:eq:dNLSE}) gives
  \[-\dot\lambda \Psi_k = \Psi_{k+1} - 2 \Psi_k + \Psi_{k-1}
  +\vert\Psi_k\vert^2 (\Psi_{k+1} + \Psi_{k-1})\]
  which is equivalent to (\ref{sr:eq:elasticCurvature}) with $\mathcal
  C = 2 - \dot\lambda$.

\end{proof}
\begin{figure}[htb]
  \begin{center}
    \centerline{\epsfxsize=0.3\hsize%
      \epsfbox{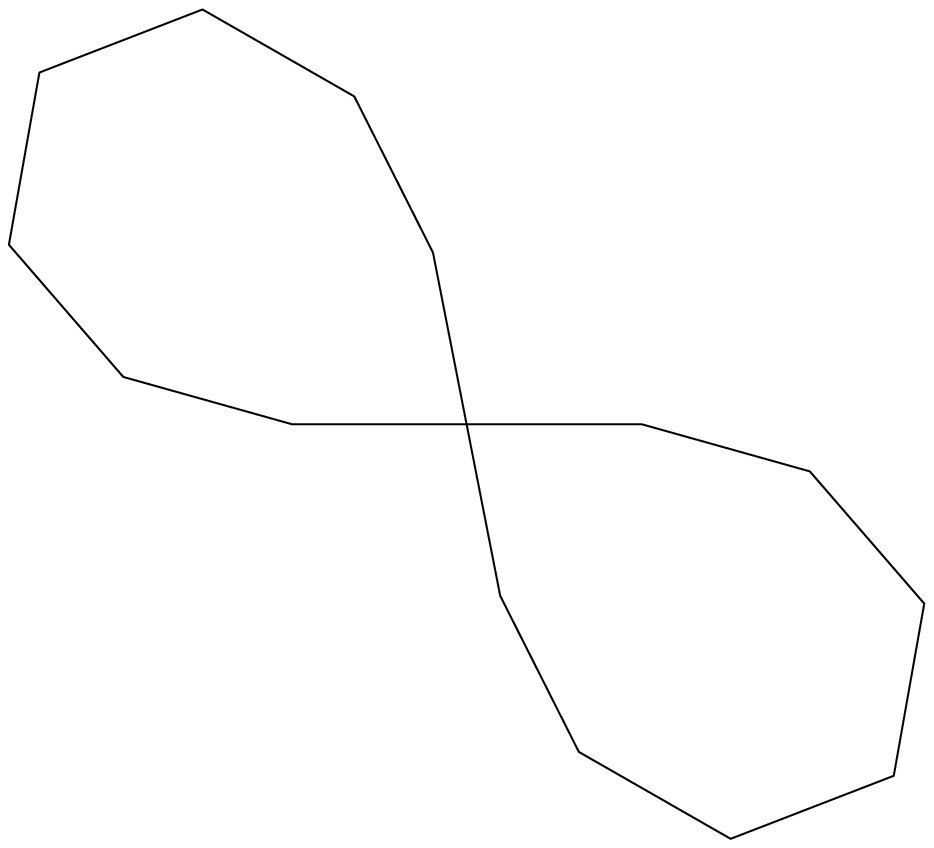}\hfil%
      \epsfxsize=0.5\hsize%
      \epsfbox{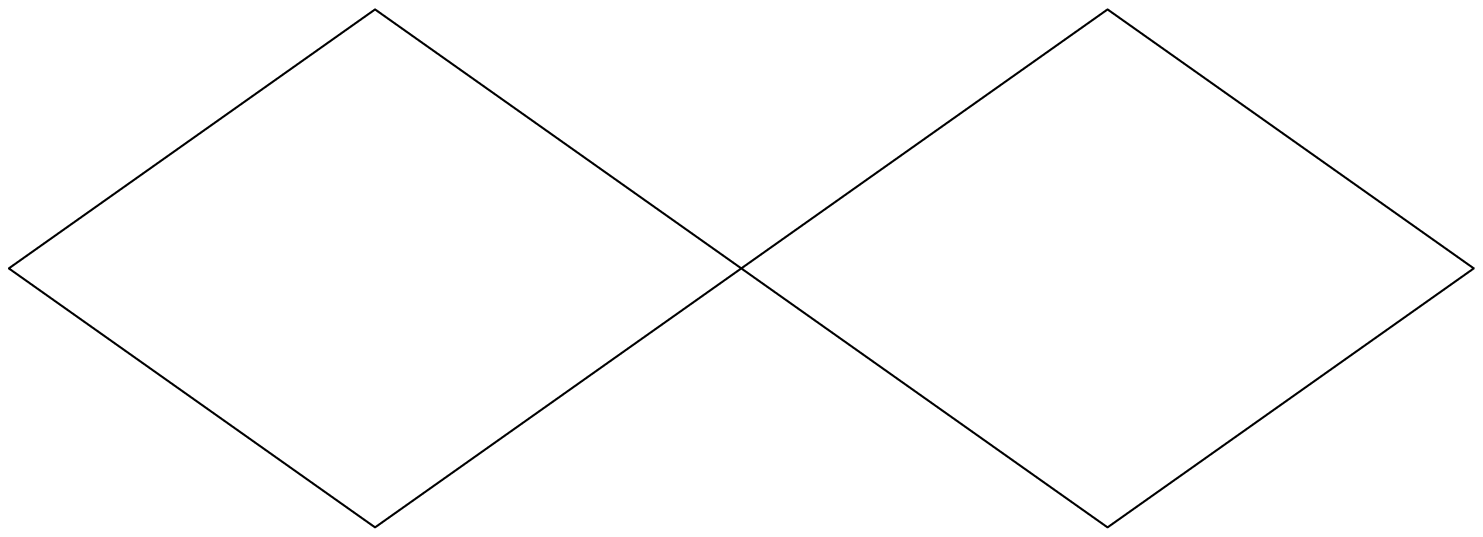}}
    \caption{Two discretizations of the elastic figure eight.}
    \label{sr:fig:Elastic8}
  \end{center}
\end{figure}
As an example Fig~\ref{sr:fig:Elastic8} shows two discretizations of
the elastic figure eight.

\subsection{B\"acklund transformations for discrete space curves and
  Hashimoto surfaces}
\label{sr:sec:discreteBaecklund}

\subsubsection{Algebraic description}

In complete analogy to Section~\ref{sr:sec:smoothAlg} we state
\begin{theorem}
\label{sr:thm:discreteDressing}
  Let $G_k$ be a solution to 
  equations (\ref{sr:eq:dIHMAuxProb})
  with $U_k$ and $V_k$ as in (\ref{sr:eq:discreteIHMLaxPair}) (i.\ e.\ $U(1)-\One$
  solves the $dIHM$ model).
  Choose $\lambda_0,s_0\in\C$. Then $ \tilde G_k(\lambda) :=
  B_k(\lambda)G_k(\lambda)$ with $B_k(\lambda) = (\One + \lambda \rho_k),
  \rho_k\in\H$ defined by the conditions that $\lambda_0,\bar\lambda_0$
  are the zeroes of $\det(B_k(\lambda))$ and 
  \begin{equation}
    \label{sr:eq:dBTcondition}
    \tilde G_k(\lambda_0){s_0\choose 1} = 0 \quad\hbox{\rm and}\quad \tilde
  G_k(\bar\lambda_0){1\choose-\bar s_0} = 0
  \end{equation}
  solves a system of the same type. In particular \[\tilde U_k(1) -\One=
  \tilde G_x(1)\tilde G^{-1}(1)-\One\] solves again the discrete
  Heisenberg magnet model 
  (\ref{sr:eq:discreteHeisenberg}).
\end{theorem}
\begin{proof}
  Analogous to the smooth case.
\end{proof}
\begin{example}
  Let us dress the (this time discrete) straight line again:
  We set $S_n \equiv \I$ and get 
  \[
  \begin{array}{rcl}
    G_n(\lambda) &=& (\One +\lambda \I)^n \exp(-2\frac{\lambda^2}{1 +
      \lambda^2}t\I)\\
    &=& 
  \quadmatrix{(1+i\lambda)^ne^{-2i\frac{\lambda^2}{1 +
  \lambda^2}t}}{0}{0}{(1 - i\lambda)^ne^{2i\frac{\lambda^2}{1 +
  \lambda^2}t}}\thep
  \end{array}\]
  After choosing $\lambda_0$ and $s_0$ and writing again 
  $\rho = \quadmatrix{a}{b}{-\bar b}{\bar a}$ we get with the
  shorthands $p = (1+i\lambda_0)^ne^{-2i\frac{{\lambda_0}^2}{1 +
  {\lambda_0}^2}t}$ and $q = (1 - i\lambda_0)^ne^{2i\frac{{\lambda_0}^2}{1 +
  {\lambda_0}^2}t}$
  \[
  \begin{array}{rcl}
    p & = & -\lambda_0(p a + s_0qb)\\
    q & = & \lambda_0(p\bar b - s_0 q \bar a)
  \end{array}\]
  which can be solved for $a$ and $b:$
  \begin{equation}
    \label{sr:eq:discreteab}
    \begin{array}{rcl}
      a & = & -\frac{\frac1{\lambda_0}\frac{\bar p}{\bar q} +
        \frac{s_0\bar s_0}{\bar \lambda_0}\frac qp}%
      {\frac{\bar p}{\bar q} + s_0 \bar s_0 \frac qp}\\[0.6cm]
      b & = & \bar s_0 \frac{\frac1{\bar\lambda_0} -
        \frac1{\lambda_0}}%
      {\frac{\bar p}{\bar q} + s_0 \bar s_0 \frac qp}\thep
    \end{array}
  \end{equation}
  Again we can write the formula for the curve $\tilde\g:$
  \[\tilde\g_n = \Im(\rho_n) + \g_n = \quadmatrix{\Im(a_n) + i n}{b_n}{-\bar
    b_n}{-\Im(a_n) - i n}\thep\]
  \begin{figure}[htbp]
    \begin{center}
      \epsfxsize=0.55\hsize\epsfbox{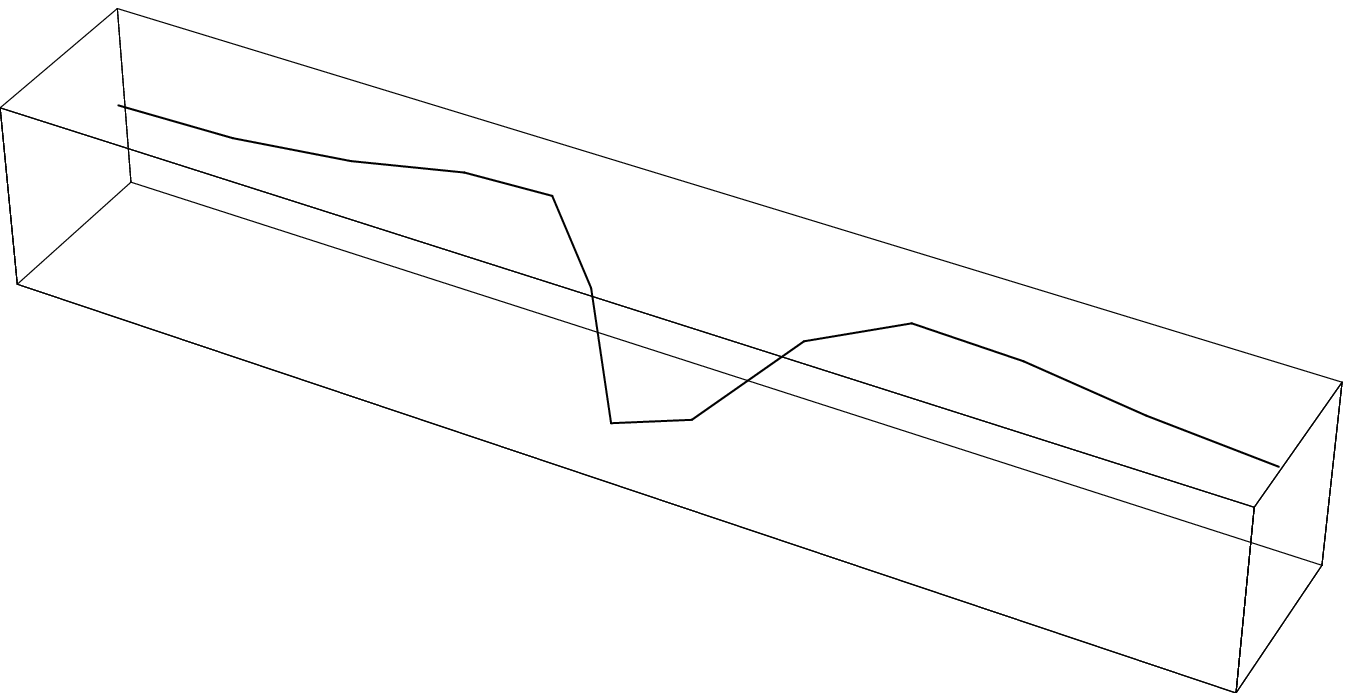}\epsfxsize=0.55\hsize\kern-0.1\hsize\raise7pt\hbox{\epsfbox{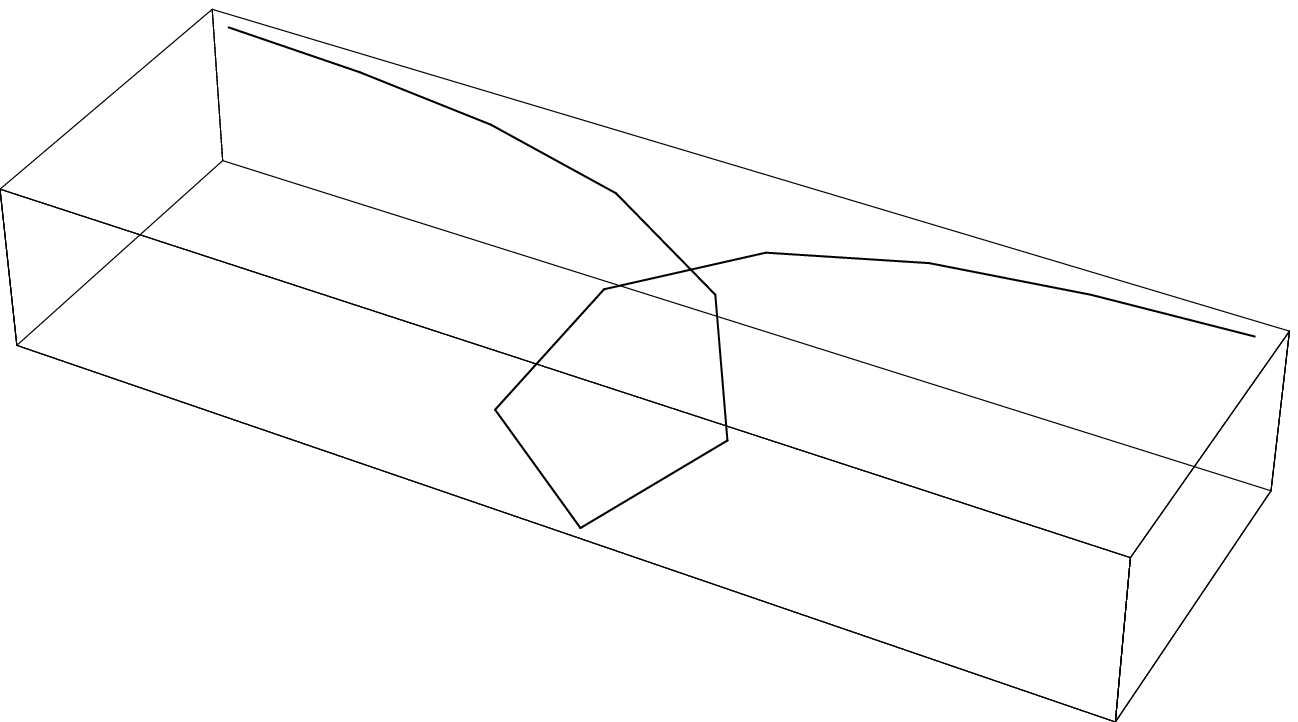}}
      
      \epsfxsize=0.55\hsize\epsfbox{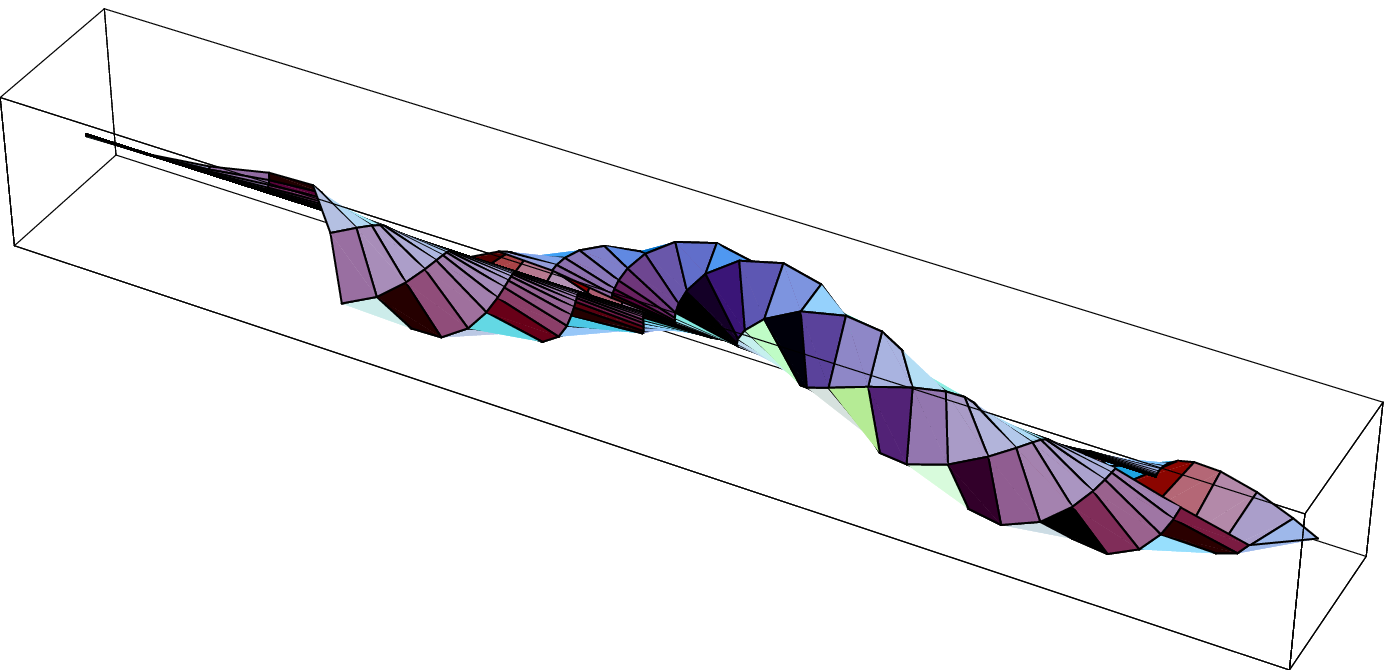}\epsfxsize=0.55\hsize\kern-0.1\hsize\raise7pt\hbox{\epsfbox{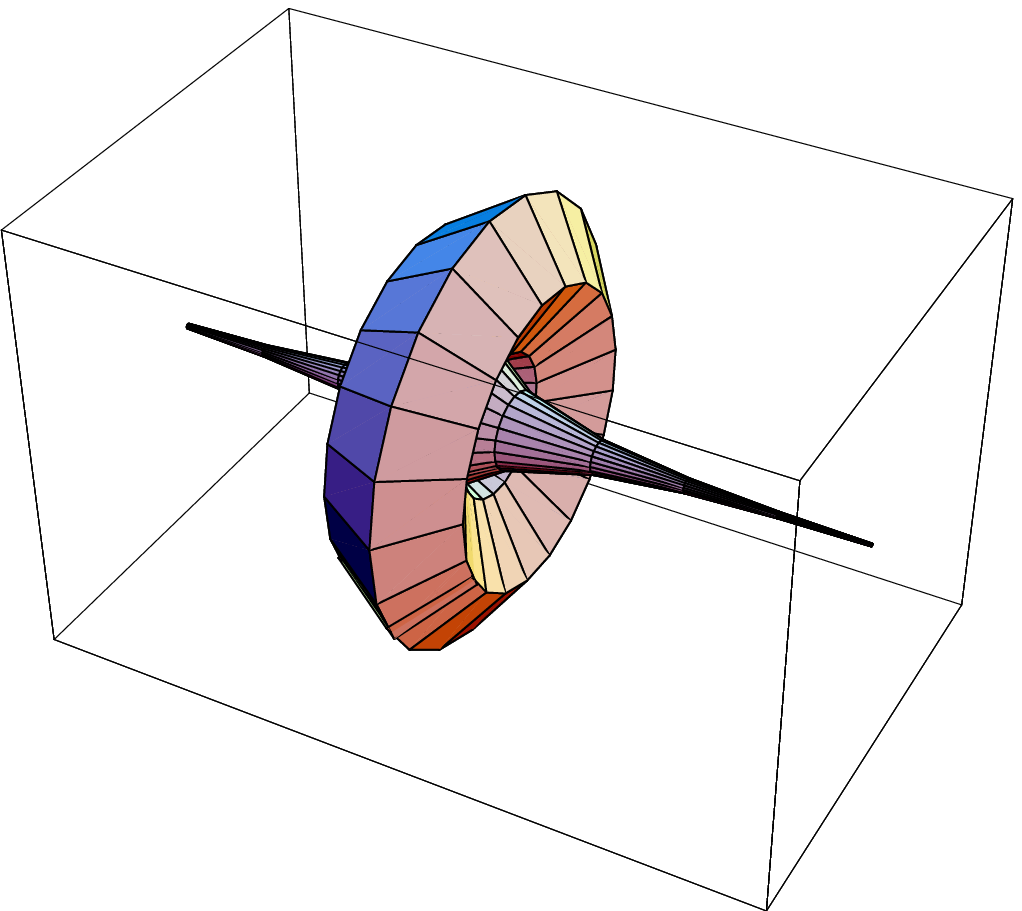}}
      \caption{Two discrete dressed straight lines and the
        corresponding Hashimoto surfaces}
      \label{sr:fig:discreteDressed}
    \end{center}
  \end{figure}
  Figure~\ref{sr:fig:discreteDressed} shows two solutions with $s_0 =
  0.5 + i$ and $\lambda_0 = 0.4 - 0.4i$ and $\lambda_0 = - 0.4i$
  respectively.  The second one is again planar. Note the strong
  similarity to the smooth examples in
  Figure~\ref{sr:fig:dressedLine}.
\end{example}
\online{
  In Figure~\ref{sr:fig:dressingApplet} you can experiment with the 
  parameters of the dressing yourself. This and the following applet
  have been programmed by Eike Preu{\ss}.
}{}
\online{
    \begin{figure}[htbp]
    \begin{center}
      \leavevmode
      \applet{third.dol.dressed.PaDressed.class}{400pt}{600pt}
      \caption{Explore the discrete dressed straight line}
      \label{sr:fig:dressingApplet}
    \end{center}
  \end{figure}
  }{}

Of course one has again a permutability theorem:

\begin{theorem}{\rm\bf (Bianchi permutability)}
  \label{sr:thm:discreteBianchi}
  Let $\tilde\g$ and $\hat\g$ be two B\"ack\-lund transforms of $\g$.
  Then there is a unique discrete Hashimoto surface $\hat{\tilde\g}$
  that is B\"acklund transform of $\tilde\g$ and $\hat\g$.
\end{theorem}
\begin{proof}
  Literally the same as for theorem~\ref{sr:thm:smoothBianchi}.
\end{proof}

\begin{figure}[htb]
  \begin{center}
    \leavevmode
    \epsfxsize=0.55\hsize\epsfbox{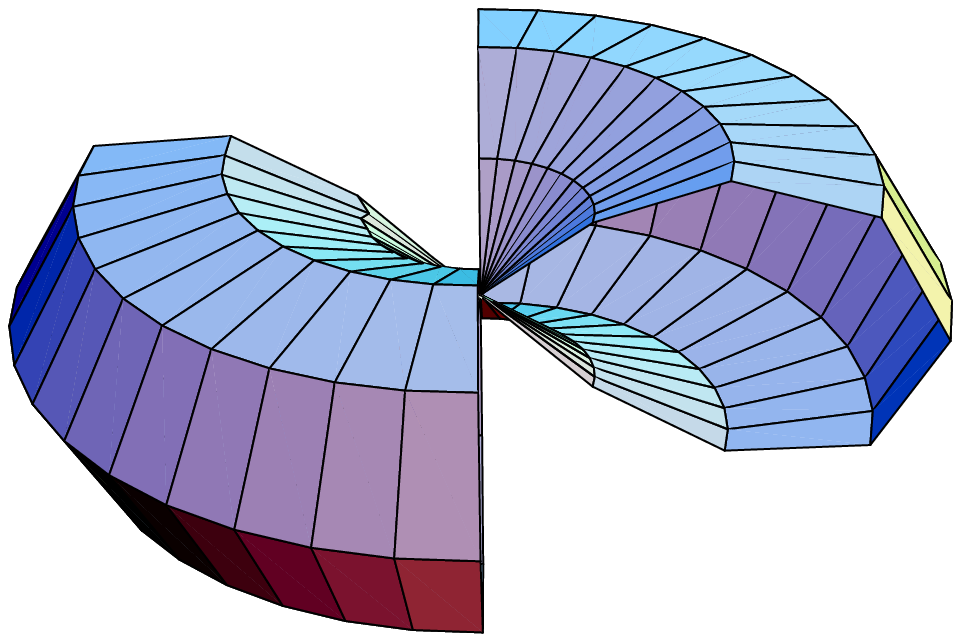}
    \caption{The Hashimoto surface from a discrete elastic eight.}
    \label{fig:sr:eightSurf}
  \end{center}
\end{figure}
\subsubsection{Geometry of the discrete B\"acklund transformation}
\label{sr:sec:DiscreteBaecklund}

In this section we want to derive the discrete B\"acklund
transformations by mimicing the twisted Traktrix construction from
Lemma~\ref{sr:contConstLengthLemma}:

Let $\g:\Z\to\Im\H$ be an discrete arclength parametrized curve. To
any initial vector $v_n$ of length $l$ there is a $S^1$-family of
vectors $v_{n+1}$ of length $l$ satisfying $\vert\g_n + v_n - (\g_{n+1}
+ v_{n+1})\vert = 1$. This is basically folding the parallelogram
spanned by $v_n$ and $S_n$ along the diagonal $S_n -
v_n$. To single out one of these new vectors let us fix the angle
$\delta_1$ between the planes spanned by $v_n$ and $S_n$ and
$v_{n+1}$ and $S_n$ (see Fig.~\ref{sr:elementaryQuad}).  This
furnishes a unique evolution of an initial $v_0$ along $\g$. The
polygon $\tilde\g_n = \g_n + v_n$ is again a discrete arclength
parametrized curve which we will call a {\em B\"acklund transform}\/ of
$\g$.

There are two cases in which the elementary quadrilaterals $(\g_n,$
$\g_{n+1},\tilde\g_{n+1}, \tilde\g_n)$ are planar. One is the
parallelogram case. The other can be viewed as a discrete version of
the Traktrix construction. 

\begin{definition}
  \label{def:sr:baecklund}
  Let $\g$ be a discrete arclength parametrized curve. Given
  $\delta_1$ and $v_0,$ $\vert v_0\vert= l$ there is a unique discrete
  arclength parametrized curve $\tilde\g_n = \g_n + v_n$ with $\vert
  v_n\vert = l$ and $\angle(\spann(v_n,S_n),\spann(v_{n+1},S_n)) =
  \delta_1$.
  
  $\tilde\g$ is called a {\em B\"acklund transform}\/ of $\g$ and
  $\hat\g = \g + \frac12 v$ is called a {\em discrete twisted
    Traktrix.}\/ for $\g$ (and $\tilde\g$).
\end{definition}
\begin{remark}
  Note that in case of $\delta = \pi$ the
  $\crt(\g,\tilde\g,\tilde\g_+,\g_+) = l^2$.
\end{remark}
Of course we will show, that this notion of B\"acklund transformation
coincides with the one from the last section.
Let us investigate this B\"acklund transformation in greater detail. 
For now we do not restrict our selves to arclength parametrized curves.
We state the following
\begin{lemma}\label{sr:moebLemma}
The map $M$ sending $v_n$ to $v_{n+1}$ in above B\"acklund
transformation is a M\"obius transformation.
\end{lemma}
\begin{proof}
Let us look at an elementary quadrilateral:
\begin{figure}
  \begin{center}
\input{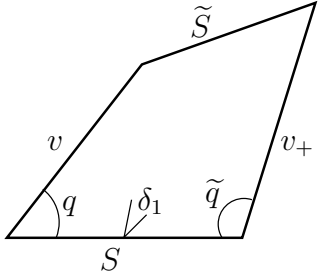}
\caption{An elementary quadrilateral of the discrete B\"acklund
transformation\label{sr:elementaryQuad}}    
  \end{center}
\end{figure}
For notational simplicity let us write $S=\g_{n+1} - \g_n,$ $\tilde S = 
\tilde\g_{n+1} - \tilde\g_n,$ $\vert S\vert = s,$ $v = v_n,$ and $v_+ = v_{n+1}$.
If we denote the angles $\angle(S,v)$ and $\angle(v_+,S)$ with $q$ 
and $\tilde q,$ we get
\begin{equation}
        e^{i \tilde q} = \frac{k e^{i q} -1}{e^{i q} - k}
        \label{sr:k-surfaceMoeb}
\end{equation}
with $k = \tan\frac{\delta_1}2\tan\frac{\delta_2}2$ and $\delta_1$
as in Fig.~\ref{sr:elementaryQuad}. $\delta_2$ is the corresponding
angle along the edge $v$.
Note that $l,$ $s,$ $k,$ $\delta_1,$ and $\delta_2$ are coupled by
\begin{equation}
        k = \tan\frac{\delta_1}2\tan\frac{\delta_2}2\quad \frac{l}{s} = 
        \frac{\sin\delta_2}{\sin\delta_1} \thep
        \label{sr:couplekld1d2}
\end{equation}
To get an equation for $v_+$ from this we need to have all vectors 
in one plane. So set $\sigma = \cos\frac{\delta_1}2 + 
\sin\frac{\delta_1}2 \frac{S}{s}$. Then conjugation with $\sigma $ is a rotation 
around $S$ with angle $\delta_1$.
If we replace $e^{i q}$ by $\frac{\sigma v\sigma^{-1}}l \left(\frac{S}{s}\right)^{-1}$
and $e^{i \tilde q}$ by 
$-\frac{S}{s}v_+^{-1}l$ equation (\ref{sr:k-surfaceMoeb}) becomes quaternionic 
but stays valid (one can think of it as a complex equation with different ``$i$'').
Equation (\ref{sr:k-surfaceMoeb}) now reads
\[
\frac{v_+ S}{ls} = \frac{\frac{s}{l} \sigma v\sigma^{-1} S^{-1} 
-k}{\frac{ks}l \sigma v\sigma^{-1} S^{-1} - 1}\thep
\]

We can write this in homogenous coordinates: $\H^2$ carries a natural right $\H$-modul
structure, so one can identify a point in $\HP^1$ with a quaternionic line in $\H^2$
by $p \cong (r,s) \iff p = rs^{-1}$. In this picture our equation gets
\[
\left(\begin{array}{c} \frac1{ls} v_+ S\\1\end{array}\right) \lambda= 
\quadmatrix{\frac{s}l \sigma}{-k S\sigma}{\frac{ks}l\sigma}{-S\sigma} 
\left(\begin{array}{c} v \\1\end{array}\right)\thep
\]
Bringing $ls$ and $S$ on the right hand side gives us finally the matrix
\begin{equation}
        \mathcal M := \quadmatrix{\frac1{k}\sigma}{-\frac{l}sS\sigma}{\frac1{ls} S\sigma}{\frac1k\sigma}\thep
        \label{sr:quaternionicMoeb}
\end{equation}
Since we know that this map sends a sphere of radius $l$ onto itself, we 
can project this sphere stereographically to get a complex matrix. The
matrix 
\[
P = \quadmatrix{2 \I}{-2l}{\frac{\J}l}{\K}
\]
projects $l S^2$ onto $\C$. Its inverse is given by
\[
P^{-1} = -\frac14 \quadmatrix{\I}{2 l \J}{\frac1l}{2\K}\thep
\]
One easily computes
\begin{equation}
        \mathcal M_\C = P\mathcal M P^{-1} = -\frac14\quadmatrix{\nu + 
        i\Re(S\I)}{2 l \Im(S\I)\J}{-\frac1{2l} \bar{\Im(S\I)\J}}{\nu -i\Re(S\I)}
        \label{sr:complexMoeb}
\end{equation}
with $\nu = i s\frac{\tan\frac{\delta_1}2 i 
-\frac1k}{\frac{\tan\frac{\delta_1}2}k i -1}$. This completes our proof.
\end{proof}
\begin{remark}
\begin{list}{--}{\labelsep1.5em}
\item Using equation 
(\ref{sr:couplekld1d2}) one can compute 
\begin{equation}
\label{sr:lambdaFormula}
\nu = 
s\tan\frac{\delta_1}2\frac{1-k^2}{\tan^2\frac{\delta_1}2 +k^2} + i l =
l\tan\frac{\delta_2}2\frac{1-k^2}{\tan^2\frac{\delta_2}2 +k^2} + i l.
\end{equation}
So the real part of $\nu$ is invariant under the change 
$s\leftrightarrow l,$ $\delta_1\leftrightarrow \delta_2$.
Therefore instead of thinking of $\tilde S$ as an transform of $S$ with parameter 
$\nu$ one could view $v_+$ a transform of $v$ with parameter $\nu 
+ i(s-l)$.
\item One can gauge $\mathcal M_\C$ to get rid of 
the off-diagonal $2l$ factors
\[ M = \quadmatrix{\frac1{\sqrt{2l}}}{0}{0}{\sqrt{2l}}\mathcal 
M_\C\quadmatrix{\sqrt{2l}}{0}{0}{\frac1{\sqrt{2l}}}\thep\]
Then we can write in abuse of notation
\begin{equation}
        M = \nu \One - S
        \label{sr:MafterGauge}
\end{equation}
Here $\nu \One$ is no quaternion if $\nu$ is complex. The 
eigenvalues of $\mathcal M_\C$ and $M$ clearly coincide and $M$
obviously coincides with the Lax matrix $U_k$ of the dIHM model in equation
(\ref{sr:eq:discreteIHMLaxPair}) up to a factor $\frac1\nu$ with
$\lambda = -\frac1\nu$.
\end{list}
\end{remark}


As prommised the next lemma shows that the geometric B\"acklund
transformation discussed in this section coincides with the one from
the algebraic description.
\begin{lemma}
\label{sr:theCrutialLemma}
Let $S,v\in\Im\H$ be nonzero vectors , $\vert v\vert = l$,
$\tilde S$ and $v_+$ be the evolved vectors in the sense of our
B\"acklund transformation with parameter 
$\nu$ ($\Im\nu = l$). then 
\begin{equation}
        (\lambda\One + \tilde S)(\lambda\One + \Re\nu +v) = (\lambda\One + 
        \Re\nu + v_+)(\lambda\One + S)
        \label{sr:zeroCurvatureIHS}
\end{equation}
holds for all $\lambda$.
\end{lemma}
\begin{proof}
Comparing the orders in $\lambda$ on both sides in equation 
(\ref{sr:zeroCurvatureIHS}) gives two equations
\begin{eqnarray}
        \tilde S + \Re\nu +v & = & \Re\nu + v_+ + S
        \label{sr:zcFirstOrder}  \\
        \tilde S (\Re\nu + v) & = & (\Re\nu + v_+)S\thep
        \label{sr:zcSecondOrder}
\end{eqnarray}
The first holds trivially from construction the second gives
\[\Re\nu  = (v_+S - \tilde Sv)(\tilde S - S)^{-1}\thep\]
This can be checked by elementary calculations using equation 
(\ref{sr:lambdaFormula}) for the real part of $\nu$.
\end{proof}

Like in the continuous case we can deduce that $\Im(\rho_n)= v_n =
\tilde\g_n - \g_n$ which gives us the constant distance between the
original curve $\g_n$ and its B\"acklund transform $\tilde\g_n$.


\section{The doubly discrete Hashimoto flow}
\label{sr:sec:ddHashimoto}
From now on let $\g:\Z\to\Im\H$ be periodic or have at least periodic 
tangents $S_n = \g_{n+1} - \g_n$ with period $N$ (we will see later that
rapidly decreasing boundary conditions are valid also).
As before let $\tilde\g$ be a B\"acklund transform of $\g$ with initial 
point $\tilde\g_0 = \g_0 + v_0,$ $\vert v_0\vert = l$.
As we have seen the map sending $v_n$ to 
$v_{n+1}$ is a M\"obius transformation and therefore the map sending 
$v_0$ to $v_N$ is one too. As such it has in general two but at least one 
fix point. Thus starting with one of them as initial point the
B\"acklund transform 
$\tilde\g$ is periodic too (or has periodic tangents $S$). Clearly this 
can be iterated to get a discrete evolution of our discrete curve $\g$.

\begin{lemma}
  Let $\g$ be a discrete curve with periodic tangents $S$ of period
  $N$. Then the tangents $\tilde S$ of a dressed curve $\tilde \g$
  with the parameters $\lambda_0$ and $s_0$ are again periodic if and
  only if the vector $(1,s_0)$ is an eigenvector of the monodromy
  matrix $G_N(\lambda)$ at $\lambda = \lambda_0$.
\end{lemma}
\begin{proof}
  We use the notation from Theorem~\ref{sr:thm:discreteDressing}.
  Since $\tilde \g_n - \g_n = v_n = \Im(\rho_n)$ and since $B(\lambda)
  = \One +\lambda \rho$ is completely determined by $\lambda_0$ and
  $v$ we have, that $B_0(\lambda) = B_N(\lambda)$. On the other hand
  on can determine $B(\lambda)$ by $\lambda_0$ and $s_0$.  Since
  $G_0(\lambda) = \One$ condition~\ref{sr:eq:dBTcondition} says that
  $1\choose s_0$ and $G_n(\lambda_0){1\choose s_0}$ must lie in $\ker
  B_0(\lambda_0)$.
\end{proof}

A Lax representation for this evolution is given by equation
(\ref{sr:zcSecondOrder}) which is basically the Bianchi permutability
of the B\"acklund transformation.

In the following we will show that for the special choice $l=1$ and
$\delta_1 \approx \frac\pi2$ the resulting evolution can be viewed as
a discrete smoke ring flow. More precisely one has to apply the
transformation twice: once with $\delta_1$ and once with $-\delta_1$.
In \cite{HO99b} it is shown, that under this evolution the
complex curvature of the discrete curve solves the doubly discrete
NLSE introduced by Ablowitz and Ladik \cite{AL76b}, which of course is
an other good argument.

\begin{proposition}
A M\"obius transformation that sends a disc into its inner has a fix point 
in it.
\end{proposition}
\begin{proof}
  For the M\"obius transformation $M$ look at the vector field $f$
  given by $f(x) = M9x) - x$. This must have a zero.
\end{proof}
\begin{figure}
  \begin{center}
    \input{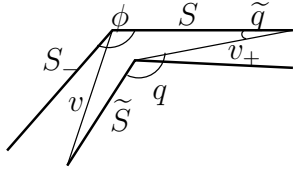}
    \caption{An elementary quadrilateral if $l = 1$ and $\delta_1\approx 
      \frac\pi2$\label{sr:elementaryQad2}}
  \end{center}
\end{figure}
Now we show the following
\begin{lemma}
If $\angle(-S_-,v)\leq \epsilon,$ $\epsilon$ sufficiently small, there
exists a $\delta_1$ such that $\angle(-S,v_+)<\epsilon$.
\end{lemma}
\begin{proof}
With notations as in Fig~\ref{sr:elementaryQad2}
we know $e^{i\tilde q} = \frac{k e^{iq} -1}{k -e^{i q}}$ and $q\in[\phi 
-\epsilon,\phi+\epsilon]$ giving us
\[
2 i \sin\tilde q = 2 i\Im e^{i\tilde q} = 2 i\frac{(k^2 - 
1)\sin(\phi\pm\epsilon)}{(k^2+1)-2k\cos(\phi\pm\epsilon)}
\]
which proofs the claim since $k$ goes to 1 if $\delta_1$ tends to $\frac\pi2$.
\end{proof}

Knowing this one can see that an initial $v_0$ with 
$\angle(-S_{N-1},v_0)\leq\epsilon$ is mapped to a $v_N$ with 
$\angle(-S_{N-1},v_N)<\epsilon$. Above Proposition gives that there must 
be a fix point $p_0$ with $\angle(-S_{N-1},p_0)<\epsilon$.

\begin{figure}[htbp]
  \begin{center}
    \centerline{\epsfxsize=0.5\hsize\epsfbox{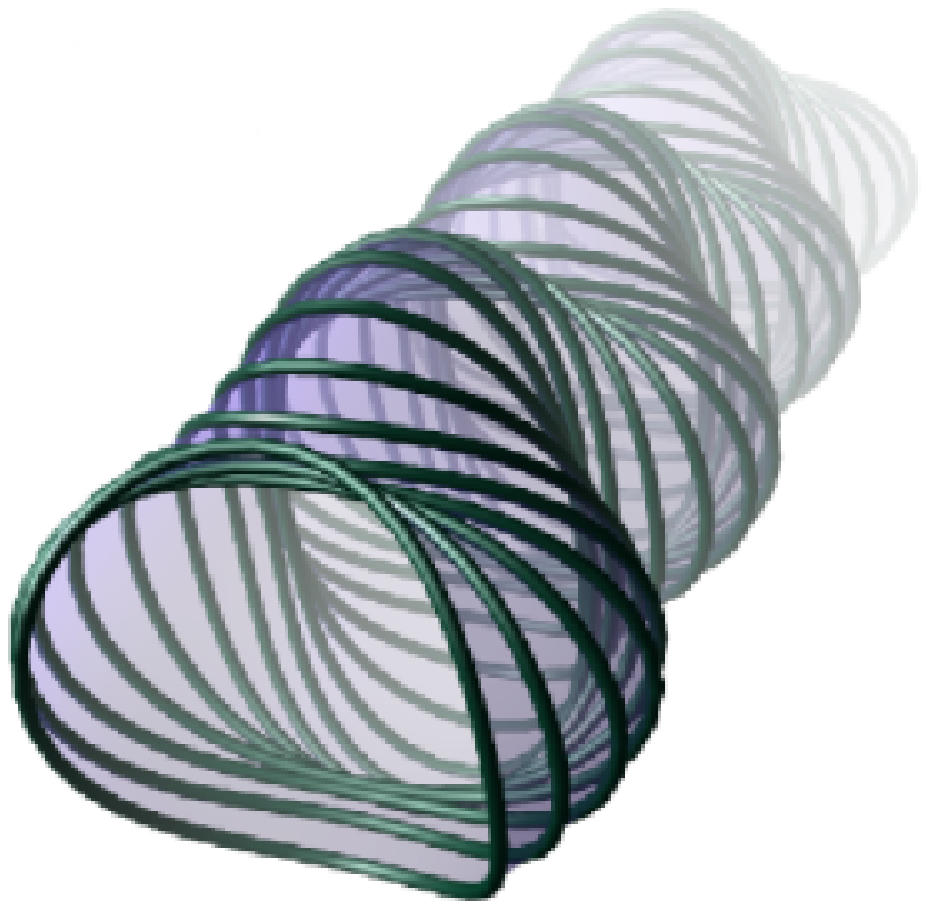}\hfil%
      \epsfxsize=0.45\hsize\epsfbox{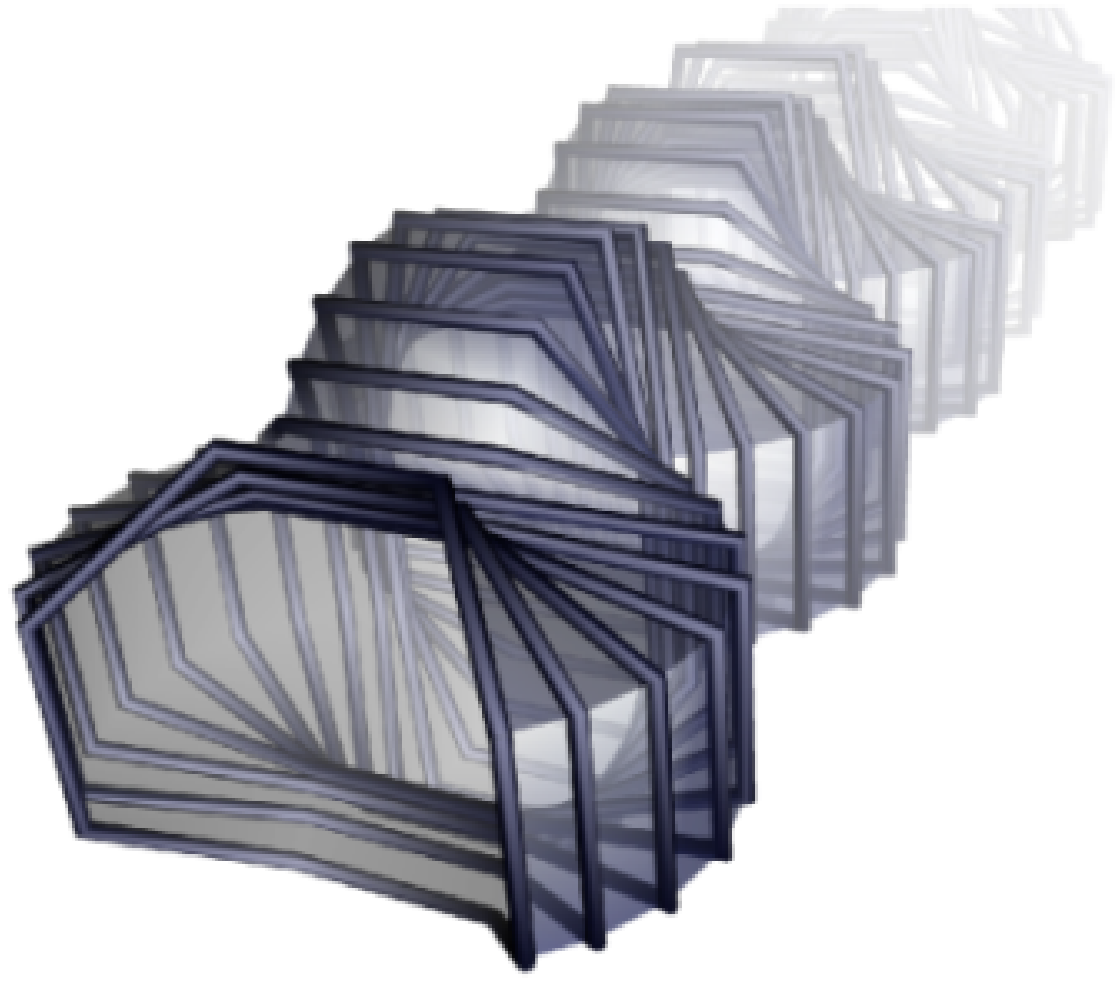}}
    \caption{An oval curve under the Hashimoto flow and the discrete evolution of its discrete pendant.}
    \label{sr:fig:flowGeneral}
  \end{center}
\end{figure}

But if $p_n\approx -S_{n-1}$ we get $\tilde\g_n\approx \g_{n-1}$ and 
$\tilde\g_n - \g_{n-1}$ is close to be orthogonal to 
$\spann(S_{n-2},S_{n-1})$. 
So it is a discrete version of an evolution in binormal direction 
---plus a shift. To get rid of this shift, one has to do the
transformation twice but with opposite sign for $\delta_1$.
Figure~\ref{sr:fig:flowGeneral} shows some stages of the smooth
Hashimoto flow for an oval curve and the discrete evolution of its
discrete counterpart\online{ and in Figure~\ref{fig:sr:evolApplet} you
  can experiment with the doubly discrete evolution of a planar curve:
  Just drag the points of the polygon as you like and then switch to
  step evolution or view the resulting surface.}{.} In general the double B\"acklund transformation
can be viewed as a discrete version of a linear combination of
Hashimoto and tangential flow---this is emphasized by the fact that
the curves that evolve under such a linear combination by rigid motion
only coincide in the smooth and discrete case:

\online{
\begin{figure}[htb]
  \begin{center}
    \leavevmode
    \applet{third.dol.flow.PaFlow.class}{400pt}{600pt}
    \caption{Discrete Hashimoto flow and surface from planar polygons}
    \label{fig:sr:evolApplet}
  \end{center}
\end{figure}
}{}



\subsection{discrete Elastic Curves}
\label{sr:elasticSection}
As a spin off of the last section one can easily show, that the
elastic curves defined in Section~\ref{sr:sec:discreteHashimoto} as
curves that evolve under the Hashimoto flow by rigid motion only do
the same for the doubly discrete Hashimoto flow. Again we will use the
evolution of the complex curvature of the discrete curve. We mentioned
before that in the doubly discrete case the complex curvature evolves
with the doubly discrete NLSE given by Ablowitz and Ladik
\cite{AL76b,HO99b}.

We start by quoting a special case of their result which can be
summarized in the following form (see also \cite{SU97})

\begin{theorem}[Ablowitz and Ladik]
\label{sr:AblowitzLadik}
  given 
  \[L_n(\mu) = \quadmatrix{\mu}{q_n}{-\bar q_n}{\mu^{-1}}\]
  and $V_n(\mu)$ with the following $\mu$--dependency:
  \[V_n(\mu) = \mu^{-2} {V_{-2}}_n + \mu^{-1} {V_{-1}}_n +
  {V_{0}}_n + \mu^{1} {V_{1}}_n + \mu^{2} {V_{2}}_n\thep\]
Then the zero curvature condition $V_{n+1}(\mu)L_n(\mu) = \tilde
L_n(\mu) V_n(\mu)$ gives the following equations:
\begin{equation}
  \label{sr:alevolution}
  \begin{array}{rcl}
    (\tilde q_n - q_n)/i &=& \alpha_+ q_{n+1} - \alpha_0 q_n +
    \bar\alpha_0 \tilde q_n - \bar\alpha_+ \tilde q_{n-1}\\
    &&+ (\alpha_+ q_n
    \mathcal A_{n+1} - \bar\alpha_+ \tilde q_n \bar{\mathcal A}_n) \\
    &&+ (-\bar\alpha_- \tilde q_{n+1} +
    \alpha_- q_{n-1})(1 + |\tilde q_n|^2 )\Lambda_n\\
    \mathcal A_{n+1} - \mathcal A_n&=& \tilde q_n \bar{\tilde q}_{n-1}
    - q_{n+1}\bar q_n\\
    \Lambda_{n+1}(1+|q_n|^2) &=& \Lambda_{n}(1+|\tilde q_n|^2)
  \end{array}
\end{equation}
with constants $\alpha_+,$ $\alpha_0$ and $\alpha_-$.
\end{theorem}
In the case of
periodic or rapidly decreasing boundary conditions the 
natural conditions $ \mathcal A_n \rightarrow 0,$ and
$\Lambda_n\rightarrow 1$ 
for $n\rightarrow \pm\infty$ give formulas for $\mathcal A_n$ and
$\Lambda_n$:
\[ \mathcal A_n = q_n \bar q_{n-1} + \sum_{j = j_0}^{n-1} (q_j\bar
q_{j-1} - \tilde q_j\bar{\tilde q}_{j-1})\]
\[ \Lambda_n = \prod_{j = j_0}^{n-1}\frac{1 + |\tilde q_j|^2}{1+|q_j|^2}\]
with $j_0 = 0$ in the periodic case and $j_0 = -\infty$ in case of
rapidly decreasing boundary conditions.

\begin{theorem}
  The discrete elastic curves evolve by rigid motion under the doubly
  discrete Hashimoto flow.
\end{theorem}
\begin{proof}
  Evolving by rigid motion means for the complex curvature of a
  discrete curve, that it must stay constant up to a possible global phase,
  i.\ e.\ $\tilde\psi_n = e^{2  i \theta}\Psi_n$. Due
  to Theorem~\ref{sr:AblowitzLadik} the evolution equation for
  $\psi_n$ reads
  \[
  \begin{array}{rcl}
    \frac{(\tilde \Psi_n - \Psi_n)}i &=& \alpha_+ \Psi_{n+1} - \alpha_0 \Psi_n +
    \bar\alpha_0 \tilde \Psi_n - \bar\alpha_+ \tilde \Psi_{n-1} +
    (\alpha_+ \Psi_n 
    \mathcal A_{n+1} \\&&- \bar\alpha_+ \tilde \Psi_n \bar{\mathcal A}_n) 
    + (-\bar\alpha_- \tilde \Psi_{n+1} +
    \alpha_- \Psi_{n-1})(1 + |\tilde \Psi_n|^2 )\Lambda_n\\
  \end{array}
  \]
  Using $e^{-i\theta}\tilde\psi_n = e^{i\theta}\Psi_n$ gives $\Delta_n = 1,$
  $\mathcal A_n = e^{2 i\theta}\Psi_n\bar\Psi_{n-1},$ and finally
  \[
    2\left(\sin\theta +\Re(e^{i\theta}\alpha_0)\right)\frac{\Psi_n}{1 +
    \vert\Psi_n\vert^2} =\] \[= \left(e^{i\theta}\alpha_+ +
    \bar{e^{i\theta}\alpha_-}\right)\Psi_{n+1} + \left
    ( \bar{e^{i\theta}\alpha_+} + e^{i\theta}\alpha_- \right)\Psi_{n-1}\thep
  \]
  So the complex curvature of curves that move by rigid motion solve 
  \begin{equation}
    \label{sr:eq:ddElasticRod}
    \mathcal C \frac{\Psi_n}{1 +
      \vert\Psi_n\vert^2} = e^{i \mu} \Psi_{n+1} +e^{-i \mu}\Psi_{n-1}
  \end{equation}
  with some real parameters $\mathcal C$ and $\mu$ which clearly
  holds for discrete elastic curves.
\end{proof}
\begin{remark}
  The additional parameter $\mu$ in eqn (\ref{sr:eq:ddElasticRod})
  is due to the fact that the Ablowitz Ladik system is the general
  double B\"acklund transformation and not only the one with
  parameters $\nu$ and $-\bar\nu$. This is compensated by the extra
  torsion $\mu$ and the resulting curve is in the associated family of
  an elastic curve. These curves are called {\em elastic rods}
  \cite{BS98}.
\end{remark}

\subsection{B\"acklund transformations for the doubly discrete \\Hashimoto
  surfaces}
\label{sr:sec:ddDressing}

Since the doubly discrete Hashimoto surfaces are build from B\"ack\-lund
transformations themselves the Bianchi permutability theorem
(Theorem~\ref{sr:thm:discreteBianchi}) ensures that the B\"acklund
transformations for discrete curves give rise to B\"acklund
transformations for the doubly discrete Hashimoto surfaces too. Thus
every thing said in section~\ref{sr:sec:discreteBaecklund} holds in
the doubly discrete case too.

\subsection*{Conclusion}
We presented an integrable doubly discrete Hashimoto or Heisenberg
flow, that arises from the B\"acklund transformation of the (singlely)
discrete flow and showed how the equivalence of the discrete and
doubly discrete Heisenberg magnet model with the discrete and
doubly discrete nonlinear Schr\"odinger equation can be understood
from the geometric point of view. The fact that the stationary
solutions of the dNLSE and the ddNLSE coincide stresses the strong
similarity of the both and the power of the concept of integrable
discrete geometry. 

Let us end by giving some more figures of examples of
the doubly discrete Hashimoto flow.


\begin{figure}[htbp]
 \begin{center}
   \epsfxsize=\hsize\epsfbox{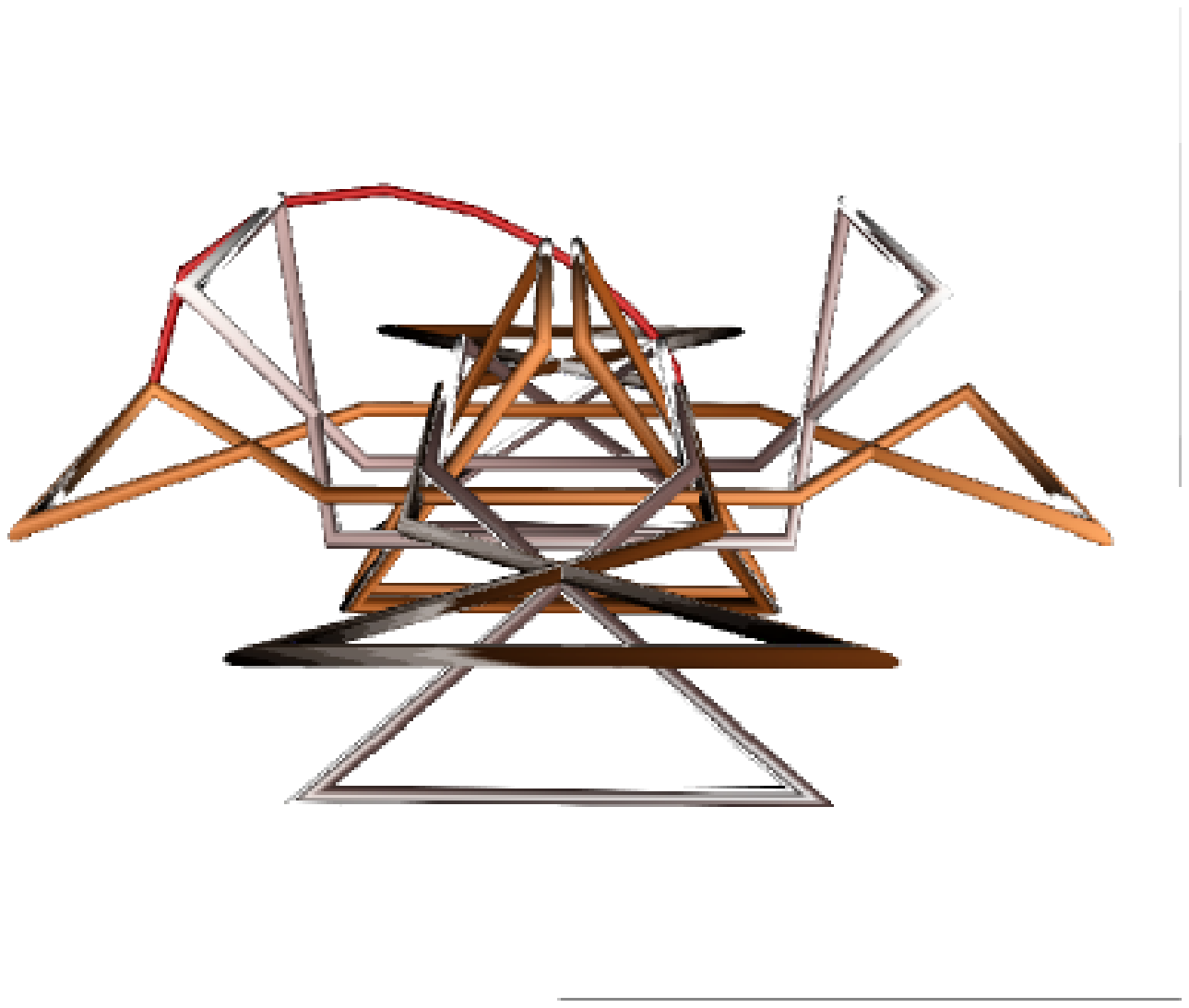}
   \caption{A discrete double eight that gives a Hashimoto torus. The red line is the trace of one vertex.}
   \label{sr:fig:doubleEight}
 \end{center}
\end{figure}
\begin{figure}[htbp]
 \begin{center}
   \epsfxsize=\hsize\epsfbox{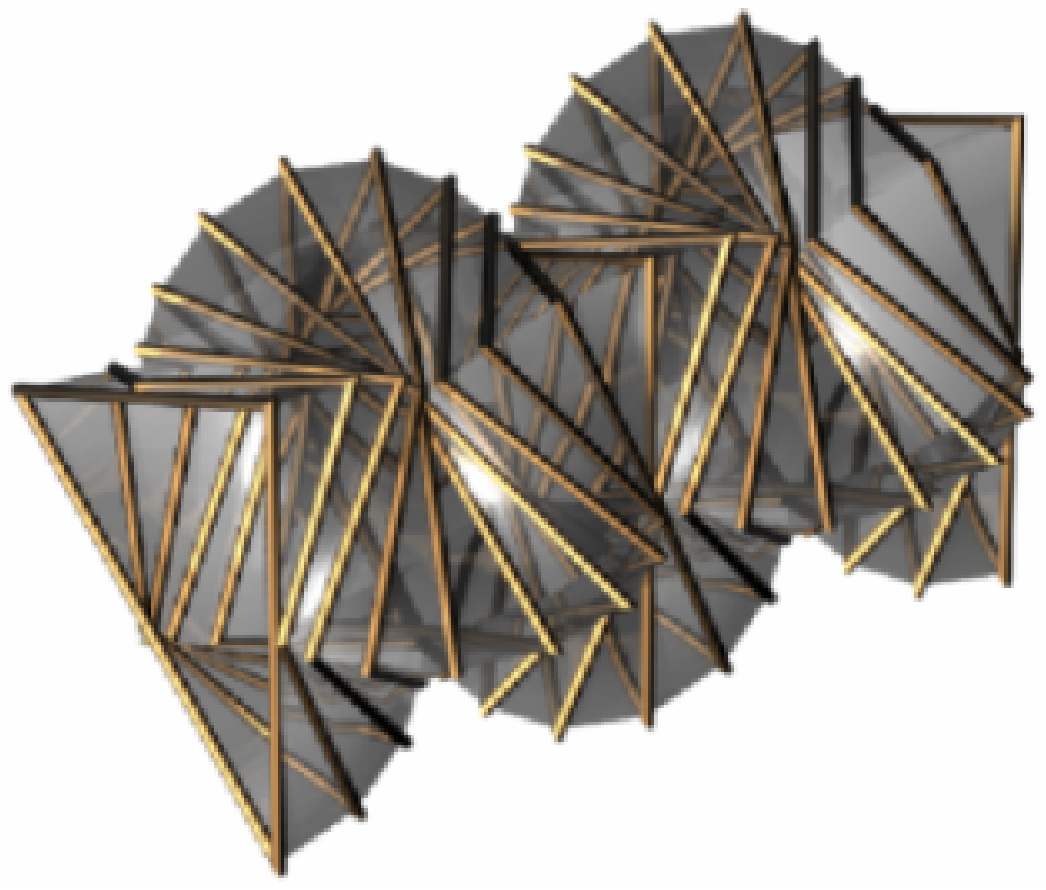}
   \caption{The doubly discrete Hashimoto flow on a equal sided triangle with subdivided edges.}
   \label{sr:fig:kaleidoZyklus}
 \end{center}
\end{figure}


\bibliographystyle{plain}
\bibliography{discrete}

\bigskip

\noindent{Tim Hoffmann\\
  \small Fachbereich Mathematik MA 8-5\\[-3pt]
  \small Technische Universit\"at Berlin\\[-3pt]
  \small Stra\ss e des 17. Juni 136\\[-3pt]
  \small10623 Berlin\\[-0pt]
  \small email: timh@sfb288.math.tu-berlin.de\\[-3pt]
  \small phone: +49-30-314 25784\\[-3pt]
  \small fax: +49 - 30 - 314 21577 }
\end{document}